\documentclass[preprint,12pt]{elsarticle}

\usepackage[utf8]{inputenc}
\usepackage{ragged2e}
\usepackage{bbm}
\usepackage{amsmath}
\usepackage{graphicx}
\usepackage{verbatim}
\usepackage{booktabs}
\usepackage{amssymb}
\usepackage{float}
\usepackage{xcolor}
\usepackage[utf8]{inputenc}
\usepackage{amssymb}
\usepackage{appendix}
\usepackage[ruled,vlined]{algorithm2e}
\usepackage{algorithmicx}
\usepackage{setspace}
\usepackage{hyperref}

\usepackage{lineno}

\newcommand{\R}{\mathbb R}

\newcommand{\bsY}{\boldsymbol{Y}}
\newcommand{\bsy}{\boldsymbol{y}}
\newcommand{\bsX}{\boldsymbol{X}}
\newcommand{\bsx}{\boldsymbol{x}}

\newcommand{\bsV}{\boldsymbol{V}}
\newcommand{\bsU}{\boldsymbol{U}}
\newcommand{\bsu}{\boldsymbol{u}}

\newcommand{\bsz}{\boldsymbol{z}}
\newcommand{\bsf}{\boldsymbol{f}}

\newcommand{\bso}{\boldsymbol{\omega}}
\newcommand{\bsO}{\boldsymbol{\Omega}}

\newcommand{\bsF}{\boldsymbol{F}}

\newcommand{\bstheta}{\boldsymbol{\theta}}

\newcount\Comments  
\newcommand{\kibitz}[2]{\ifnum\Comments=1\textcolor{#1}{#2}\fi}

\begin{document}
\begin{frontmatter}

\title{ State Space Kriging model for emulating complex nonlinear dynamical systems under stochastic excitation }

\author{Kai Cheng$^{a*}$,  Iason Papaioannou$^{a}$, MengZe Lyu$^{b}$ and Daniel Straub$^{a}$ }

\affiliation{organization={Engineering Risk Analysis Group, Technical University of Munich},
 addressline={Theresienstr. 90}, 
city={ Munich},
 postcode={80333}, 
 country={Germany}}
\affiliation{organization={College of Civil Engineering, Tongji University},
 addressline={Siping Rd. 1239}, 
city={Shanghai},
 postcode={200092}, 
 country={ China}}
 
\begin{abstract}

We present a new surrogate model for emulating the behavior of complex nonlinear dynamical systems with external stochastic excitation. The model represents the system dynamics in state space form through a sparse Kriging model. The resulting surrogate model is termed  state space Kriging (S2K) model. Sparsity in the Kriging model is achieved by  selecting an informative training subset from the observed time histories of the state vector and its derivative with respect to time. We propose a tailored technique for designing the training time histories of state vector and its derivative, aimed at enhancing the robustness of the S2K prediction. We validate the performance of the S2K model with various benchmarks. The results show that S2K yields accurate prediction of complex nonlinear dynamical systems under stochastic excitation with only a few training time histories of state vector.
\end{abstract}


\end{frontmatter}


\section{Introduction}

Dynamical systems are widely used in modern engineering and applied science for modeling complex underlying physical phenomena \cite{benner2015survey}. With the increase of computational power, numerical simulation offers a feasible way to study and predict the behavior of complex dynamical systems. However, the response of complex dynamical systems is governed by uncertainties, due to the stochastic external excitations, uncertain boundary conditions, and natural variability of system properties \cite{li2009stochastic}. To obtain effective prediction,  these uncertainties must be accounted for. To this end, uncertainty quantification of stochastic dynamical systems has gained particular interest in the last few decades \cite{mai2016surrogate,lyu2024decoupled}.

In the context of uncertainty quantification, sampling-based simulation methods, i.e., Monte Carlo simulation (MCS) and various variance reduction methods \cite{kanjilal2021cross,bayer1999importance,au2007application}, are generally used to propagate uncertainties from system inputs to system response quantities of interest. These methods are robust, but they are infeasible when only a small number of expensive simulations are affordable or available. To address this issue, surrogate models have been widely used to construct computationally efficient approximations of the expensive computational model. Various surrogate modelling techniques have developed, including: Gaussian process regression (aka, Kriging) \cite{oliver1990kriging,rasmussen2003gaussian,kleijnen2009kriging}, support vector regression (SVR) \cite{smola2004tutorial,drucker1996support,cheng2021adaptive}, polynomial chaos expansion (PCE) \cite{xiu2002wiener,blatman2011adaptive,cheng2018adaptive,papaioannou2019pls}, neural networks \cite{abiodun2018state,raissi2019physics,shi2024novel,hao2024multi}. These surrogate modelling techniques are powerful for approximating the behavior of traditional static \textquotedblleft black-box\textquotedblright models, but exhibit difficulties when applied to dynamical systems under stochastic excitation \cite{mai2016surrogate}. In these problems, the number of input parameters due to discretization of the stochastic external excitations, often represented in terms of a white noise process, can be extremely high. Most common surrogate models suffer from the \textquotedblleft curse of dimensionality\textquotedblright, since effective learning typically requires at least twice as many samples as the number of input parameters.
To address this issue, a standard technique is to insert existing surrogate modelling techniques into the framework of nonlinear auto-regressive with exogenous input (NARX) modelling \cite{menezes2008long,mai2016surrogate,bhattacharyya2020kriging}. The NARX model is a powerful system identification technique, which is established based on the principle of causality, i.e., it assumes that the system response quantity of interest at the current time instant is only affected by its previous several response values and the  current and past several values of external excitations. Based on this cause-consequence effect, the current response quantity of interest is assumed to be a function of the response values at past multiple time instants and the input excitation at the current and previous multiple instants. This function can be emulated through application of  existing surrogate modelling techniques.  Although the NARX model has proved its effectiveness in several structural dynamical problems \cite{mai2016surrogate,wan2023feature,worden2018confidence,bhattacharyya2020kriging}, it is only an empirical model, and there is no rigorous way to select the model hyper-parameters, e.g., the time lags of both input and output. In addition, the NARX model struggles to emulate the response of complex highly nonlinear dynamical problems. Recently, a manifold NARX (mNARX) model \cite{schar2024emulating}  has been proposed to address the above problems, in which a transformation function is adopted to map the input into a problem-aware manifold, which is expected to be more suitable for constructing the NARX model than in original input space. However, physical information of the dynamical system under investigation is required to construct such input manifold in mNARX. If no physical information is available, the mNARX degenerates to the traditional NARX. 
 
In the present work, we propose to emulate complex nonlinear dynamical systems under stochastic excitation in their state space form with the Kriging model, termed as the state space Kriging (S2K) model. The state space representation of a stochastic dynamical system can be considered as a multi-input multi-output (MIMO) function, in which the input is the state vector and the external excitation, and the output is the derivative of the state vector. The Kriging model is utilized to learn every component of this MIMO function separately. To overcome the inefficiency of the Kriging surrogate for large training data set, we introduce an active learning algorithm to select a reduced training set from the training time histories, resulting in a sparse Kriging model. By learning the state space form of a dynamical system, the S2K model avoids the curse of dimensionality 
resulting from the discretization of the stochastic external excitation. Numerical examples demonstrate that the S2K model is effective for emulating various complex nonlinear dynamical systems with stochastic excitation with only a few training time histories of the state quantities, and outperforms the NARX model.

The layout of this paper is as follows. In section 2, we review the fundamental definitions of nonlinear dynamical systems under stochastic excitation and the basic idea of the NARX model. The S2K model is presented in Section 3, together with the active learning algorithm and the technique for designing the training time history of state vector. In Section 4, several benchmarks are used to assess the performance of our method. The paper concludes with final remarks in Section 5.

\section{Background}
\label{sec::Background}

A dynamical system is generally represented by a high-order differential equation. It can be equivalently transformed into multiple coupled first-order equation systems by introducing new state variables, a form known as the state space representation. 
In the present work, we consider a general nonlinear dynamical system under stochastic excitation in its state space form \cite{brunton2016discovering}, which can be generally expressed as
 \begin{eqnarray}
\dot\bsX(t)  = \bsf(\bsX(t),\bsU(t)), \, \mathrm{with}\enspace \bsX(0) = \bsx_0,
 \label{eq:dynamical_system}
 \end{eqnarray}
 where $\bsX(t)=[X_1(t),...,X_n(t)]^{\rm T} \in \mathbbm{R}^{n} $ is the state vector at time $t$; $\bsx_0$ is the initial condition; $\dot\bsX(t)$ is the derivative of $\bsX(t)$ with respect to $t$; $\bsU(t)=[U_1(t),...,U_m(t)]^{\rm T}\in \mathbbm{R}^{m}$ is the external stochastic excitation vector acting on the structure; $\bsf(\cdot)$ is the $n$-dimensional nonlinear vector function . 
 
 For every realization of the time history $\bsu(t)$ of $\bsU(t)$, the time history of the state vector $\bsx(t)$ of Eq. \eqref{eq:dynamical_system} exists and uniquely depends on the initial conditions $\bsx_0$, and one can apply various numerical discretization methods, e.g., the Runge-Kutta method, to find the solutions $\bsx(t_i)(i=1,...,N_t)$ at $N_t$ time steps within the time period of interest $[0,T]$. In general, $\bsf(\cdot)$ is a computational expensive  "black-box" model, which makes UQ of complex  dynamical systems under stochastic excitation a computationally demanding task.  

  The NARX model is a popular method for UQ of dynamical system under stochastic excitation \cite{mai2016surrogate,schar2024emulating,wan2023feature}. It is a system identification technique developed based on the discrete-time step representation of the dynamical system. Given a discrete time history of the state vector $\bsx(t_1),...,\bsx(t_{N_t})$ corresponding to a discrete realization of the stochastic excitation $\bsU(t)$, namely, $\bsu(t_1),...,\bsu(t_{N_t})$, the NARX model represents the response quantity of interest  $x(t_i)\in \mathbbm{R}$  at the current time instant as a function of its past values and the input excitation values at the current and previous instants as
\begin{eqnarray}
x(t_i)= g(\bsu(t_i),\bsu(t_{i-1}),...,\bsu(t_{i-n_u}),x(t_{i-1}),...,x(t_{i-n_x})) + \epsilon_t,
\label{eq:narx}
\end{eqnarray}
where $n_u$ and $n_x$ represent the maximum excitation and response time lags; $\epsilon_t$ is the residual of the NARX model; $g(\cdot)$ is the underlying model to be learned. For dynamical systems, the function $g(\cdot)$ is usually learned by polynomial model or Kriging model \cite{worden2018confidence}.  

The NARX model can be interpreted as an empirical model that tries to capture the behaviour of the stochastic dynamical system with a low-dimensional auto-regressive surrogate model. Its accuracy highly depends on the choice of the time lags $n_u$ and $n_x$, and there is no rigorous way to determine them. For problems with long memory, both $n_u$ and $n_x$ become large, and the NARX model is high-dimensional. Moreover, NARX gives poor prediction for strongly nonlinear dynamical systems, which we also confirm in Section \ref{subsec:bouc}.

\section{Methodology}

 In this section, we introduce the proposed algorithm for learning  dynamical systems under stochastic excitation in state space form. We also present the active learning algorithm for selecting the informative training set and the technique for designing the training time history of state vector and its derivative. 

\subsection{Kriging model for learning dynamical systems }

By denoting $\bsY(t)=\dot\bsX(t)$, the original dynamical system in Eq. \eqref{eq:dynamical_system} can be expressed as
 \begin{eqnarray}
 \bsY(t)  = \bsf(\bsX(t),\bsU(t))
\label{eq:transformed_dynamical_system}.
 \end{eqnarray}

    Note that $\bsf(\cdot): \mathbbm{R}^{m+n} \rightarrow \mathbbm{R}^{n} $ is a deterministic function, which maps the state vector $\bsX(t)  \in \mathbbm{R}^{n} $ and the external excitation $\bsU(t)  \in \mathbbm{R}^{m}$ to $\bsY(t) \in \mathbbm{R}^{n}$ at time instant $t$. The input dimension of $\bsf(\cdot)$ depends on the dimension of the state vector $n$ of the associated dynamical system and the cardinality of the external excitation vector $m$. For every realization of the time history $\bsu(t)$ of the excitation $\bsU(t)$, the time history of the corresponding state vector $\bsx(t)$ and its derivative  $\bsy(t)$ over the time period of interest $[0,T]$ can be estimated numerically. 

In this work, we employ the Kriging model \cite{kleijnen2009kriging,williams1995gaussian} to approximate the state space representation $\bsf(\cdot)$. Note that other existing surrogate modelling techniques, such as PCE, SVR and neural network could also be implemented in this context. 
As in the ordinary Kriging model, it is assumed that every component $y_i(t)$ of $\bsy(t)$ is a  Gaussian process with constant mean $\beta_i$:
\begin{eqnarray}
y_i(t) = \beta_i + \sigma_i^2 Z(\bso(t)), \,\ \mathrm{for}\enspace i =1,...,n,
\end{eqnarray}
where $\bso(t)=[\bsx(t),\bsu(t)]^{\rm
T}\in \mathbbm{R}^{n+m}$, $\sigma_i^2$ is the process variance, and $Z(\bso(t))$ is a stationary Gaussian process with zero mean and unit variance. The correlation function of $Z(\bso(t))$ is given by
$R(\bso(t),\bso(t');\bstheta_i)$, which describes the correlation between $\bso(t)$ and $\bso(t')$ at time $t$ and $t'$, with the correlation lengths in various coordinate directions being controlled by the hyper-parameter vector $\bstheta_i\in \mathbbm{R}^{m+n}$. Moreover, we assume that $y_i(t)$ is independent of $y_j(t)$ for all $i \neq j$, which allows learning $\bsy(t)$ component-wisely.
   
Given the experimental design $\bsO_t=[\bso(t_1),...,\bso(t_N)]\in \mathbbm{R}^{(n+m)\times N}$, the joint distribution $y_i(t)$ and  $\bsY_i=[y_i(t_1),...,y_i(t_N)]^{\rm
T} \in \mathbbm{R}^{N}$ corresponding to $\bsO_t$ is multivariate Gaussian given by
\begin{equation*}
 \left[\begin{array}{cc} y_i(t) \\ \bsY_i  
 \end{array}\right] \sim \mathcal {N}\left(\left[\begin{array}{cc} \beta_i  \\ \beta_i \boldsymbol{F}\end{array}\right],\sigma_i^2\left[\begin{array}{cc} 1 & \boldsymbol{r}^{\rm T}(\bso(t)) \\ \boldsymbol{r}(\bso(t)) & \boldsymbol{R}\end{array} \right] \right),
\end{equation*}
where $\bsF=[1,...,1]^{\rm T} \in \mathbbm{R}^{N}$ and
\begin{eqnarray*}
\begin{aligned}
\ \boldsymbol{r}(\bso(t))&:= [R(\bso(t_{i_1}),\bso(t);\bstheta_i)]_{i_1} 
\in \mathbbm{R}^{N}, \\
\boldsymbol{R}&:=[R(\bso(t_{i_1}),\bso(t_{i_2});\bstheta_i)]_{{i_1},{i_2}}
\in \mathbbm{R}^{N\times N}.
\end{aligned}
\end{eqnarray*}

The predictive distribution of $y_i(t)$ given the observation $\bsY_i=\bsy_i$ is still Gaussian, i.e.,
\begin{eqnarray}
\hat y_i(t)  \sim \mathcal {N}\left(\mu_i(t), s_i^2(t)\right),
 \label{eq:conditional_predictor}
 \end{eqnarray}
where the predictive mean and predictive variance are given by

\begin{equation}
       \mu_i(t) = \beta_i +r^\mathrm{T}(\bso(t))\boldsymbol{R}^{-1}(\bsy_i-\beta_i\bsF)  ,
       \label{eq:mean}
     \end{equation}
     and
     \begin{equation}
     s_i^2(t) = \sigma_i^2 \left(1-r^\mathrm{T}(\bso(t))\boldsymbol{R}^{-1}r(\bso(t))+\left(1-\bsF^\mathrm{T}\boldsymbol{R}^{-1}r(\bso(t))\right)^2/(\bsF^\mathrm{T}\boldsymbol{R}^{-1}\bsF)\right).
     \label{eq:variance}
    \end{equation}

The predictive mean is the kriging surrogate model prediction, and the predictive variance is used for measuring the predictive uncertainty. 
The prediction $\hat\bsy(t)$ of the vector $\bsy(t)$ is obtained by stacking $\hat y_i(t)(i=1,...,n)$, namely, $\hat\bsy(t)=[\hat y_1(t),...,\hat y_n(t)]^{\rm T}$.  $\hat\bsy(t)$ is a Gaussian vector, with mean $\boldsymbol{m}(t) = [m_1(t),..., m_n(t)]^{\rm T}$ and diagonal covariance matrix with diagonal vector $\boldsymbol{s}^2(t) = [s_1^2(t),...,s_n^2(t)]^{\rm T} $, namely
\begin{eqnarray*}
\hat\bsy(t) = \left\{\begin{array}{ccc} \hat y_1(t)\\ \vdots \\ \hat y_n(t) \end{array} \right\} \sim \mathcal {N}\left(\left\{\begin{array}{ccc} \mu_1(t)\\ \vdots \\ \mu_n(t) \end{array} \right\}, \left[\begin{array}{ccc} s_1^2(t) & \cdots & 0  \\ \vdots  & \ddots &  \vdots   \\ 0 & \cdots &s_n^2(t) \end{array} \right] \right).
 \end{eqnarray*}

In the current work, the $Mat\acute{e}rn$-$5/2$ correlation function is used, which is given by
\begin{eqnarray*}
 R(\bso(t),\bso(t');\bstheta_i) = 
 \left(1+\sqrt{5}r +\frac{5}{3}r^2\right) {\rm exp}\left(-\sqrt{5}r\right),
\end{eqnarray*}
where $r=\sqrt{(\bso(t) -\bso(t'))^{\rm T} \boldsymbol{\Lambda}(\bso(t) -\bso(t'))}$, and $\boldsymbol{\Lambda}$ is a diagonal matrix with diagonal element $\bstheta_i$.

The hyper-parameter $\beta_i$, $\sigma_i^2$ and $\bstheta_i$ for $i$-th surrogate model are determined by the maximum likelihood estimation method \cite{williams1995gaussian}. The likelihood function of the hyper-parameters given data  reads
\begin{equation}
f(\beta_i,\sigma_i^2,\bstheta_i)=\frac{[{\rm det}\boldsymbol{R} (\bstheta_i)]^{-0.5}}{\sqrt{(2\pi\sigma_i^2)^{N}}}
{\rm exp}\left(-\frac1{2\sigma_i^2}(\bsy_i-\beta_i\bsF)^{\rm T}\boldsymbol{R}^{-1}(\bstheta_i)(\bsy_i-\beta_i\bsF)\right).
     \label{originallikelihood}
    \end{equation}
     
    After taking the logarithm of Eq. (\ref{originallikelihood}), the optimal values of $\beta_i$ and $\sigma_i^2$, depending on hyper-parameter $\bstheta$, can be derived analytically. This yields
     \begin{eqnarray}
      \beta_i(\bstheta_i) &=& (\boldsymbol{F}^\mathrm{T}\boldsymbol{R}^{-1}(\bstheta_i)\boldsymbol{F})^{-1}\boldsymbol{F}^\mathrm{T}\boldsymbol{R}^{-1}(\bstheta_i)\bsy_i,
     \label{beta} \\
     \sigma_i^2(\bstheta_i) &=& \frac1{N}(\bsy_i-\beta_i\boldsymbol{F})^{\rm T}\boldsymbol{R}^{-1}(\bstheta_i)(\bsy_i-\beta_i\boldsymbol{F}).
    \label{sigma2}
    \end{eqnarray}
    
 There is no closed-form solution for the optimal hyper-parameter $\bstheta_i$, and one has to use numerical optimization algorithms to determine its value. Substituting $\sigma_i^2(\bstheta_i)$ in Eq. \eqref{sigma2} into the likelihood function in Eq.  \eqref{originallikelihood}, we are left with minimizing the following $\bstheta_i$-dependent objective function
    \begin{equation}
     \ell(\bstheta_i)=N\ln \sigma_i^2(\bstheta_i) + \ln\left[{\rm det\boldsymbol{R}}(\bstheta_i)\right].
   \label{eq:reducedlikelihood}
   \end{equation}

 The objective function in Eq. \eqref{eq:reducedlikelihood} is generally highly nonlinear and multimodal \cite{cheng2023sliced}, and thus a global optimization algorithm is required to find the hyper-parameter. In the present work, we use the multi-starts gradient-free “Hooke $\&$ Jeeves” pattern
search method \cite{altinoz2019multiobjective} for determining the optimal hyper-parameter vector. Note that the single-start “Hooke $\&$ Jeeves” pattern
search method has been implemented in the popular DACE toolbox \cite{lophaven2002dace}, and here we adapt it to multi-starts case for improving the robustness of Kriging model. 
 
\subsection{Prediction with S2K model}
\label{sec:active_learning}
Once the S2K model is trained, it can be used to predict the time history of system response $\hat\bsx(t)$ for arbitrary realization of the time history  $\bsu(t)$ of $\bsU(t)$ by solving the following approximate dynamical system in state space form
\begin{eqnarray}
 \hat\bsy(t) = \hat\bsf(\hat\bsx(t),\bsu(t))
 \label{eq:GP_dynamic},
\end{eqnarray}
where $\hat\bsf(\cdot)$ is the Kriging approximation of $\bsf(\cdot)$.

In general, the Runge-Kutta method can be used to estimate the time history of the state vector $\hat\bsx(t_i)(i=1,...,N_t)$ at $N_t$ discretized time steps within the time period of interest $[0,T]$. 
At each time step, the Kriging model provides a  Gaussian prediction $\hat\bsy(t)$ given a realization of the current predicted state vector $\hat\bsx(t)$ and external excitation $\bsu(t)$. However, the whole time-history of $\hat\bsy(t)$ and $\hat\bsx(t)$ over the time period of interest is non-Gaussian. This is due to the recursive nature of the dynamical system, i.e., the output quantities at current time instant will be used as input to predict the output quantities at the next time step, and carry their predictive uncertainty with them. Since the surrogate model $\hat\bsf(\cdot)$ is nonlinear, the predictive distribution of the entire time history is generally 
 intractable. In practice, one can adopt MCS to approximate the predictive distribution of the state vector over the time period of interest given a realization of the excitation. This can be achieved by emulating the time history of state vector $N_{MC}$ times based on the trained S2K model, and use these time history trajectories to estimate the empirical predictive distribution.

\subsection{Sparse Kriging with active learning}

We now analyze the computational cost of the S2K model.  When training a Kriging model with $N$ training samples, one needs to determine the hyper-parameter $\bstheta_i$ by minimizing the likelihood function in Eq. \eqref{eq:reducedlikelihood}. To  numerically evaluate this function, one needs to invert the correlation matrix  $\boldsymbol{R}\in\R^{N\times N}$. When a Cholesky factorization is used to decompose $\boldsymbol{R}$, the corresponding computational costs are $\mathcal{O}(N^3)$. The computational costs for training the S2K model are therefore $\mathcal{O}(nN^3)$, which tends to be very time-consuming for large number of samples $N$ and large $n$.

 In the Kriging model, the computational cost for making a single prediction is $\mathcal{O}(N)$ for the mean in Eq. \eqref{eq:mean}, and $\mathcal{O}(N^2)$ for the variance in Eq. \eqref{eq:variance}. 
To estimate the predictive distribution of the entire time history of the state vector, one needs to predict every component of the $n$-dimensional state vector $\hat\bsx(t)$ at $N_t$ time steps,  $N_\mathrm{MC}$ times for every realization of the time history $\bsu(t)$ of $\bsU(t)$. The associated operations are $\mathcal{O}(nNN_tN_\mathrm{MC})$ for the mean and $\mathcal{O}(nN^2N_tN_\mathrm{MC})$ for the variance. In general, the sample size $N$, the time history length $N_t$ and the MCS sample size $N_\mathrm{MC}$ are all quite large, which makes it time-consuming to train the S2K model and to  estimate the predictive distribution of the entire time history of the state vector corresponding to a single realization of the time history $\bsu(t)$. 

To improve the training and prediction efficiency, we suggest constructing a sparse Kriging \cite{liu2020gaussian} model by selecting a training subset of size $N_s(N_s \ll N)$ from the available training set for every component of $\bsy(t)$. This selection is performed adaptively, through first training a Kriging model with a few initial training samples uniformly selected from the whole sample set, and then updating the Kriging model by enriching informative samples selected sequentially through maximization of the following mean square error criterion
\begin{eqnarray}
 L_i(t) = (m_i(t) - y_i(t))^2 + s_i^2(t),
\label{eq:acquistion} 
\end{eqnarray}
where the first term represents the square of the predictive bias, and the second term denotes the predictive variance. The bias term prefers points with large local predictive error, and the predictive variance term tends to select points with large global predictive uncertainty. The acquisition function in Eq. \eqref{eq:acquistion} therefore accounts for both the bias for local exploitation and the variance for global exploration. 

The proposed method selects one sample in each active learning step until the maximum of $L_i(t)/\sigma^2_y$ is less than a threshold $\delta_i$ (e.g., $10^{-5}$), where $\sigma^2_y$ is the variance of observed model response in the training set.
 The convergence threshold $\delta_i$ is a user-specified parameter, depending on the complexity of the dynamical system to be learned. With a small $\delta_i$, many samples will be selected to construct the S2K model, which improves the prediction accuracy, but decreases the prediction efficiency; with a large $\delta_i$, only a small portion of the training set will be selected to construct the S2K model, which leads to a very sparse model, but the accuracy may be low. 
Based on our numerical experiments, we suggest setting $\delta_i\in [10^{-3}, 10^{-7}]$. A higher threshold values should be chosen when the nonlinearity of $f_i(\cdot)$ is expected to be high, and vice versa. 

\subsection{Design of training time history}

The Kriging model is powerful for interpolation, but its performance degenerates when extrapolating away from the training data. To ensure the accuracy of the S2K model, the training data should cover the whole parameter space of both the state vector and the external excitation. However, for a black-box model, the parameter space of state vector is unknown beforehand, since one cannot know the minimum and maximum values of the time histories of the state vector associated with future realizations of external excitation.

To address this issue, one can collect many observed time histories of state vector corresponding to different realizations of external excitation, and select a subset of the time histories (e.g., requiring that the maximum value exceeds a pre-defined threshold) to train the surrogate model, as has been proposed for training the NARX model \cite{mai2016surrogate}. However, this is infeasible when only a few time histories of the state vector are available. In the current work, we suggest collecting the time histories of a pseudo state vector resulting from pseudo external excitation with magnified variability. By magnifying the variability of the excitation, the state vector will exhibit stronger variability, hence, it is expected that a few pseudo time histories of state vector will sufficiently populate its parameter space under the true excitation. To this end, one can magnify the standard deviation of the stochastic excitation at every time instant by a magnification factor $\sigma$ . 

\begin{figure}[htbp]
\centering
\includegraphics[width=1\textwidth, trim={10 225 10 225},clip]{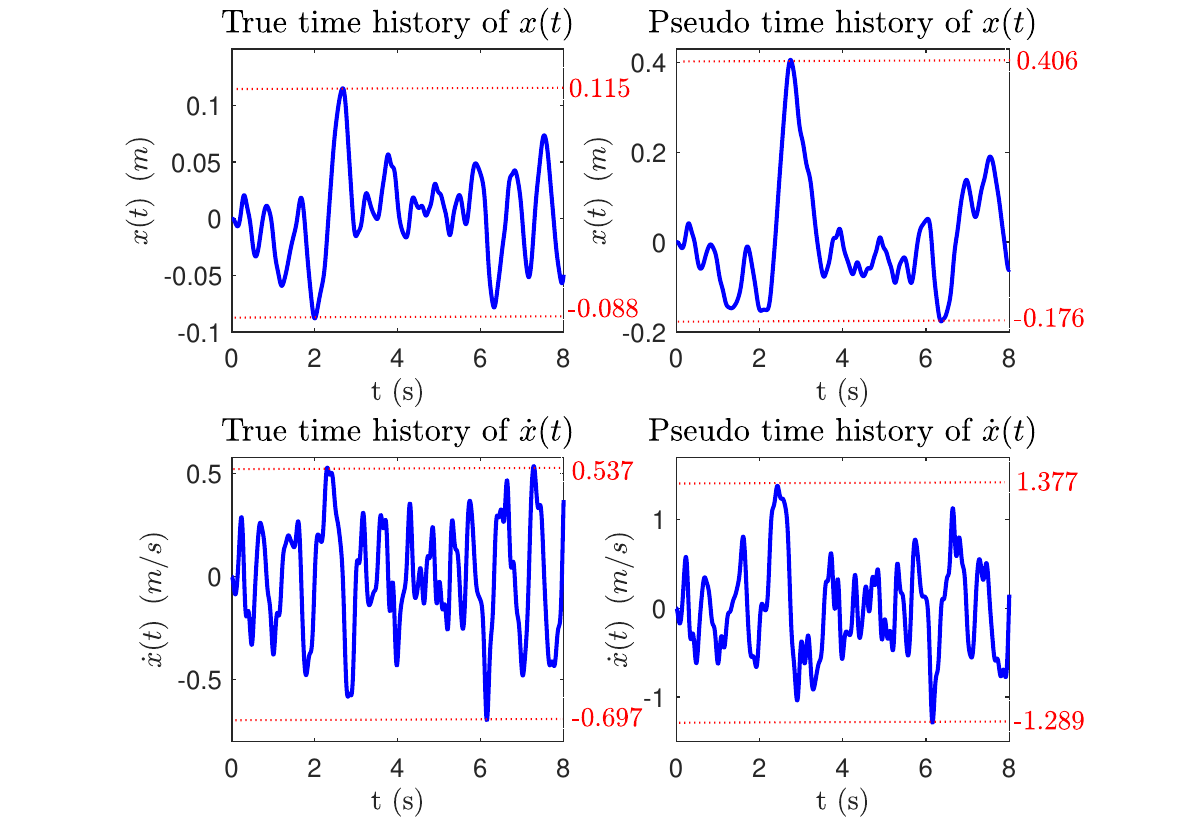}
\caption{A time history of displacement and velocity of the Bouc-Wen model with true  excitation and pseudo one with magnification factor $\sigma=2$.}  
\label{fig:design}
\end{figure}

In Fig. \ref{fig:design}, a time history of state variables (displacement and velocity) of the Bouc-Wen hysteretic oscillator (see Section \ref{subsec:bouc}) with white noise ground acceleration is presented, where the state variables are obtained with both the true excitation and a pseudo excitation with magnified variability ($\sigma=2$). It is found that the peaks and the valleys of the pseudo responses (both displacement and velocity) largely exceed the ones of the original response. Consequently, the pseudo time history is more likely to cover the whole parameter space of the state vector under the original excitation. By training the Kriging model with the pseudo training time history, extrapolation in the prediction stage can be potentially avoided, thereby improving the accuracy and robustness of the S2K model. 
Note that the magnification factor $\sigma$ is a problem-dependent parameter. A small $\sigma$ cannot guarantee full coverage of the parameter space of the state vector, while a large $\sigma$ will lead to a less accurate Kriging model, since many samples will be enriched during the active learning procedure to explore the  unnecessary region.  Numerical investigation reported in Section \ref{sec:numerical_example}  show that $\sigma \in [1.5, 2]$ is a good choice when only one training time history of the state vector is available. However, when one can afford multiple training time histories, it is beneficial to set different $\sigma $ values for each training time history. In this work, we suggest setting
\begin{equation}
\sigma_k=1+\frac{k-1}{n_t-1}, 
k=1,...,n_t
\label{eq:mixture_sigma}
\end{equation}
 where $n_t(n_t\geq 2)$ is the number training time histories. 
 In doing so, the mixture of different training time histories of state vector should better cover the entire parameter space, leading to more accurate prediction, as demonstrated in Section \ref{sec:numerical_example}. 

\subsection{Summary of the S2K algorithm}

The S2K algorithm for emulating the response of dynamical systems under stochastic excitation is summarized in Algorithm \ref{alg:algorithm3}.

\begin{small}
 \begin{algorithm}[H]
\caption{ S2K algorithm }
   \label{alg:algorithm3}
   \KwIn{  Number of training time histories $n_t$; Magnification factors $\sigma_l$ $(l=1,...,n_t)$; Convergence threshold $\delta_i(i=1,...,n)$; The state space representation $y_i(t)=f_i(\bsx(t),\bsu(t))(i=1,...,n)$; Number of initial training sample size $N_0$.  }
   \KwOut{The S2K model.}
   \begin{algorithmic}[1]
  
   \State{ Sample the training time history of excitation, and collect the training time histories of the state vector and its derivative; Assemble the entire training data set $\{\hat\bsO_t,\hat\bsy_i\}(i=1,...,n)$; Set the initial values of the learning function $L_i(t)=1(i=1,...,n)$ and $k=1$.  }\\
    \While{$k<n$}{
     {Select the $N_0$ initial training sample set $\{ \bsO_t,\bsy_k\}$  uniformly from the entire training data set $\{\hat\bsO_t,\hat\bsy_k\}$; }
     
   \While{${\rm max}(L_k(t)) >\delta_k$}{
    {2.a: Construct the Kriging model $\hat y_k(t)$ based on $\{\bsO_t,\bsy_k\}$; }
    
   {2.b: Select the most informative sample $\bsx(t^*) \in \bsO_t$ by maximizing $L_k(t)$; }
   
   {2.c: Enrich current sample set $\bsO_t=\bsO_t\cup \{\bsx(t^*),\bsu(t^*)$\} and $\bsy_k=\bsy_k \cup y_k(t^*)$;  }

    }
     { $k =k+ 1$;}
    }
   \end{algorithmic}
   \end{algorithm}
\end{small}

\section{Numerical investigations}
\label{sec:numerical_example}
We investigate the effectiveness of S2K model on four example applications. To assess the accuracy of the S2K model, the relative error  of a whole time history of the state vector corresponding to a realization of the stochastic excitation is used to define the accuracy metric as
 \begin{eqnarray}
  {\epsilon}_i =\frac{  \sum_{j=1}^{N_t}(x_i(t_j) - \hat x_i(t_j) )^2} {\sum_{j=1}^{N_t}(x_i(t_j) - \bar x_i )^2}
  \label{eq:relative_error},
  \end{eqnarray}
where $x_i(t_j)$ and $\hat x_i(t_j)$ are true and predicted $i$-th state quantity of interest at time instant $t_j$, and $\bar x_i$ is the mean over the whole time history. In addition, we use the mean value $\bar {\epsilon}_i$ of ${\epsilon}_i$ 
corresponding to $N_\mathrm{MC}$ different realizations of the excitation to measure the global accuracy of $i$-th state, namely
 \begin{eqnarray}
  \bar\epsilon_i = \frac{1}{N_\mathrm{MC}} \sum_{j=1}^{N_\mathrm{MC}}\epsilon_i^{(j)}
  \label{eq:relative_error},
  \end{eqnarray}
where $\epsilon_i^{(j)}$ is the relative error of $i$-th state quantity of interest corresponding to $j$-th realization of the stochastic excitation. We set $N_\mathrm{MC}=1000$.

For the first three example applications, we compare the performance of the S2K model with the NARX model in Eq. \eqref{eq:narx}. In the NARX model, the polynomials are used as bases, and the least angle regression technique \cite{mai2016surrogate} is used to select the most relevant basis terms. In addition, the maximal time lags of both excitation and response $n_u$ and $n_x$ are chosen equal to twice
the number of degrees of freedom of the considered dynamical system \cite{mai2016surrogate}. In these examples, both S2K and NARX model are trained 10 times independently, and the relative errors are depicted by boxplots, whereby the central mark indicates the median, and the bottom and
top edges of the box indicate the 25th and 75th percentiles, respectively. 

In both NARX and S2K model, the training time history of the state vector $\bsx(t)$ and its derivative are evaluated at equidistant time instants $t_i=i\triangle t(i=1,...,T/\triangle t)$ over the time period of interest $[0,T]$ by the Matlab solver ode89, where $\triangle t$ is the step length. Note that while equidistant discretization of the training time history is required in the NARX model, it is not essential for the S2K model. However, we apply it here for consistency.

\subsection{Quarter car model}

\begin{figure}[htbp]
\centering
\includegraphics[width=1.1\textwidth, trim={20 120 20 90},clip]{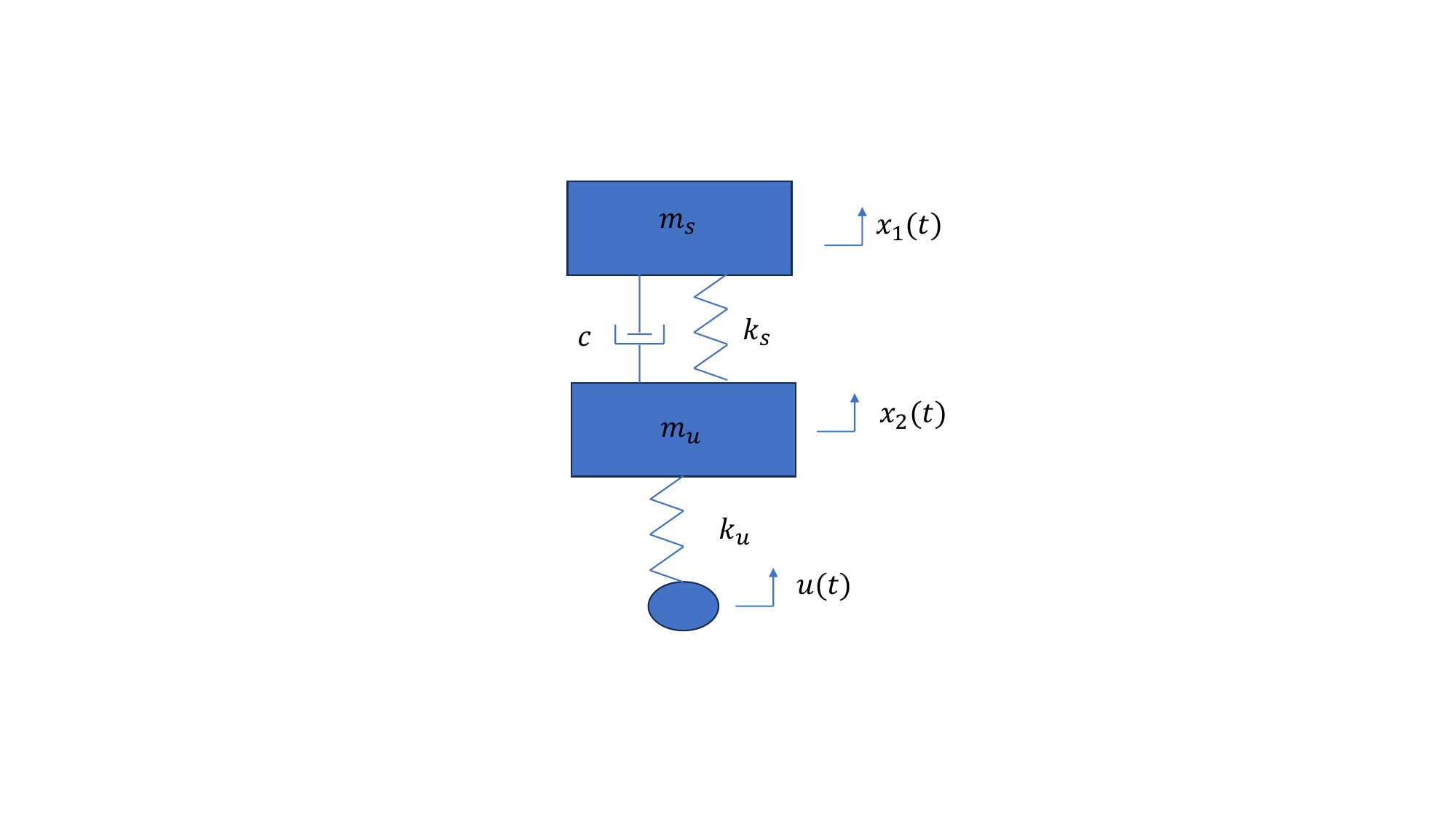}
\caption{ The quarter car model}    
\label{fig:qurater_car_model}
\end{figure}

 We first consider a quarter car model represented by a
nonlinear two degree-of-freedom system \cite{mai2016surrogate} depicted in Fig. \ref{fig:qurater_car_model}. The displacements of
the masses are governed by the following system of ordinary differential equations (ODEs), namely,
\begin{eqnarray}
\begin{aligned}
\begin{cases} 
m_s \Ddot{x}_1(t) &= -k_s(x_1(t)-x_2(t))^3 - c(\dot x_1(t)-\dot x_2(t)), \\  m_u \Ddot{x}_2(t) &= k_s(x_1(t)-x_2(t))^3 - c(\dot x_1(t)-\dot x_2(t)) + k_u(u(t)-x_2(t)),
\end{cases}
\end{aligned}
\label{eq:quarter_car}
\end{eqnarray}
where the sprung mass $m_s=22.7 \,\ \mathrm{kg}$ and the unsprung mass $m_u=42 \,\ \mathrm{kg}$ are connected by a nonlinear spring of stiffness $k_s = 1897.02 \,\ \mathrm{N/m^3}$ and a linear damper with damping coefficient $c=601.8 \,\ \mathrm{N\cdot s/m}$. An external excitation is applied to $m_u$ through a linear spring of stiffness $k_u = 1771.4 \,\ \mathrm{N/m}$. $x_1(t)$ and $x_2(t)$ are the displacements of $m_s$ and $m_u$, respectively. In this work, the excitation is modeled by a stochastic process as $u(t) = A {\rm sin}(b t)$, where the parameters $A$, $b$ follow uniform distributions, $ A \sim U(0.09,0.11) (\mathrm{m})$ and $b \sim U(1.8\pi,2.2\pi) (\mathrm{rad/s})$.

Denoting $\bsx(t)=[x_1(t),\dot x_1(t),x_2(t),\dot x_2(t)]^{\rm T}$ and $\bsy(t)=\dot\bsx(t)$, one can express the ODEs in Eq. \eqref{eq:quarter_car} in state space form as
\begin{eqnarray}
\begin{aligned}
\begin{cases} 
{y}_1(t) &  =\dot{x}_1(t), \\  {y}_2(t)  &= -\frac{k_s}{m_s}(x_1(t)-x_2(t))^3 - \frac{c}{m_s}(\dot x_1(t)-\dot x_2(t)),  \\
{y}_3(t) &  =\dot{x}_2(t), \\
{y}_4(t)  &= \frac{k_s}{m_u}(x_1(t)-x_2(t))^3 - \frac{c}{m_u}(\dot x_1(t)-\dot x_2(t)) + \frac{k_u}{m_u}(u(t)-x_2(t)). 
\end{cases}
\end{aligned}
\label{eq:quarter_car_state_sapce}
\end{eqnarray}

\begin{figure}[htbp]
\centering
\includegraphics[width=0.9\textwidth, trim={20 25 20 60},clip]{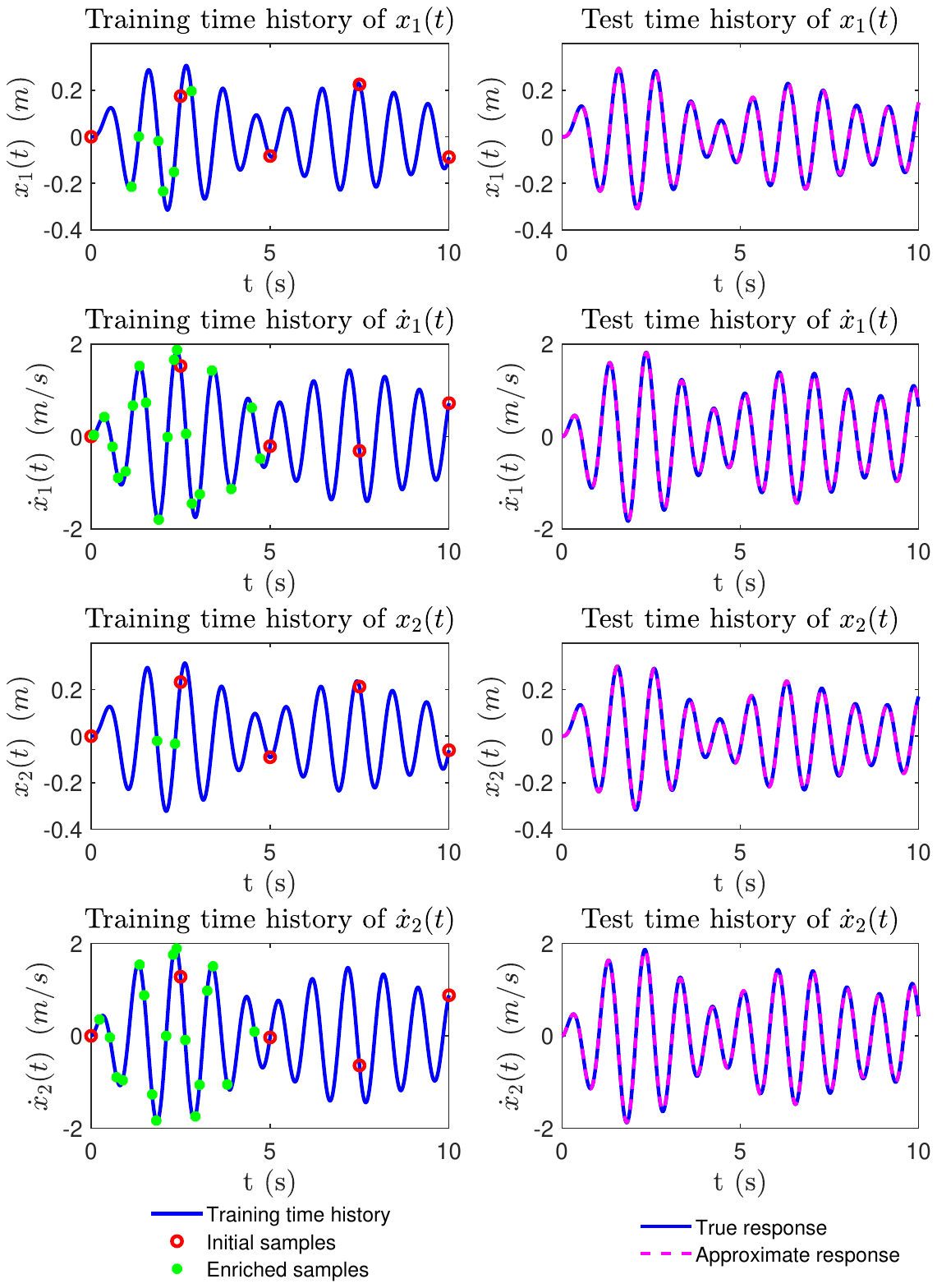}
\caption{ An example of the training and test time history of the state vector $\bsx(t)=[x_1(t),\dot x_1(t),x_2(t),\dot x_2(t)]^{\rm T}$ of the quarter car model}    
\label{fig:qurater_car_example}
\end{figure}

We use the S2K model to construct a surrogate model $ \hat\bsf(\cdot):  \mathbbm{R}^{5} \rightarrow \mathbbm{R}^{4} $ for the state space representation in Eq. \eqref{eq:quarter_car_state_sapce}, in which the input is $\bso(t) = [\bsx(t), u(t)]^{\rm T} \in \mathbbm{R}^{5}$, and the output is $\bsy(t) \in \mathbbm{R}^{4}$. 
In this example, we draw a time history of length $10$s of the excitation randomly ($A=0.0964,b=6.5248$), and the time history of the state vector and its derivative are estimated with time step $\triangle t = 0.002$s.  We use this time history of the state vector and the corresponding derivative to construct the four Kriging models, in which $5$ equidistant samples selected from the whole time history are used to train the initial Kriging models, and the active learning algorithm introduced in Section \ref{sec:active_learning} is utilized to select the informative samples sequentially with the convergence threshold $\delta=10^{-6}$. The corresponding results are depicted in Fig. \ref{fig:qurater_car_example}, where a test time history is also presented. The results demonstrate that the S2K model is highly accurate. Indeed, the relative errors of the four Kriging models corresponding to the four state quantities of the test time histories depicted in Fig. \ref{fig:qurater_car_example} are ${\epsilon}_1=6.09\times 10^ {-4}$, ${\epsilon}_2=4.54\times 10^ {-4}$, ${\epsilon}_3=6.02\times 10^ {-4}$, and ${\epsilon}_4=4.50\times 10^ {-4}$, respectively. In addition, it is observed that only a small portion of the samples are selected from the whole time history (a time history consists 5001 discrete time instants in total) after the active learning procedure. Specifically, the final sample sizes of the four Kriging models are $11$, $24$, $ 7$, and $23$, respectively. That is, the corresponding degrees of sparsity are $0.22\%$, $0.48\%$, $ 0.14\%$, and $0.46\%$.

\begin{figure}[htbp]
\centering
\includegraphics[width=0.95\textwidth,
 trim={20 260 20 280},clip]{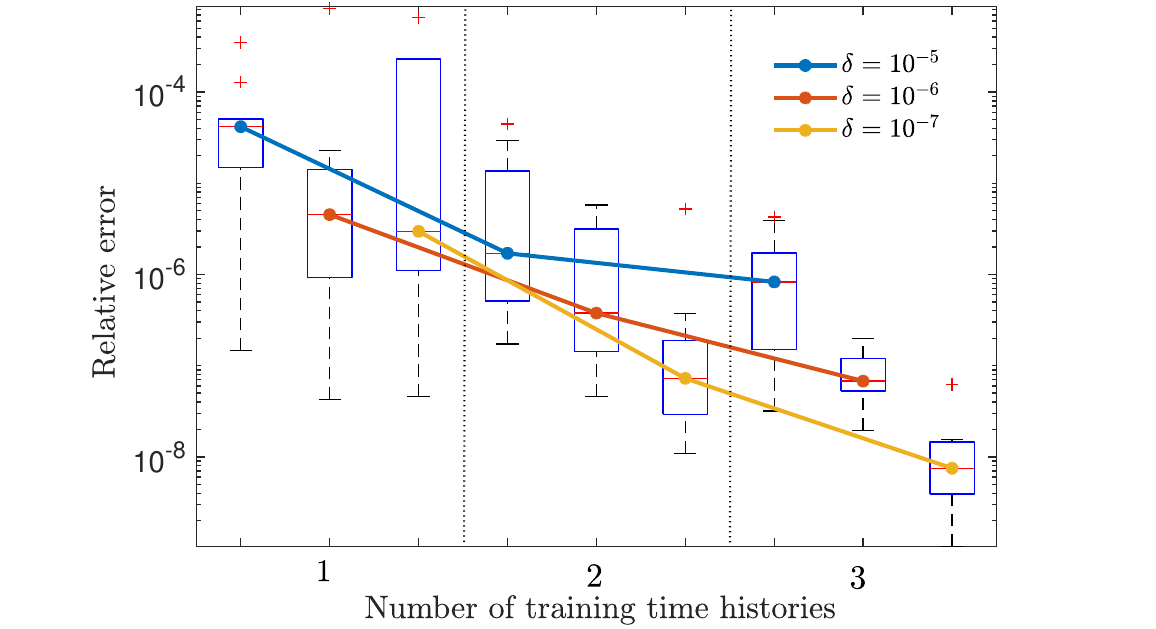}
\caption{ Boxplot of relative errors of $x_1$ (with 10 repetitions) for the quarter car model obtained with the S2K model by varying the number of training time histories with different convergence thresholds. }    
\label{fig:qurater_car_error}
\end{figure}

\begin{figure}[htbp]
\centering
\includegraphics[width=1\textwidth,trim={20 225 20 230},clip]
{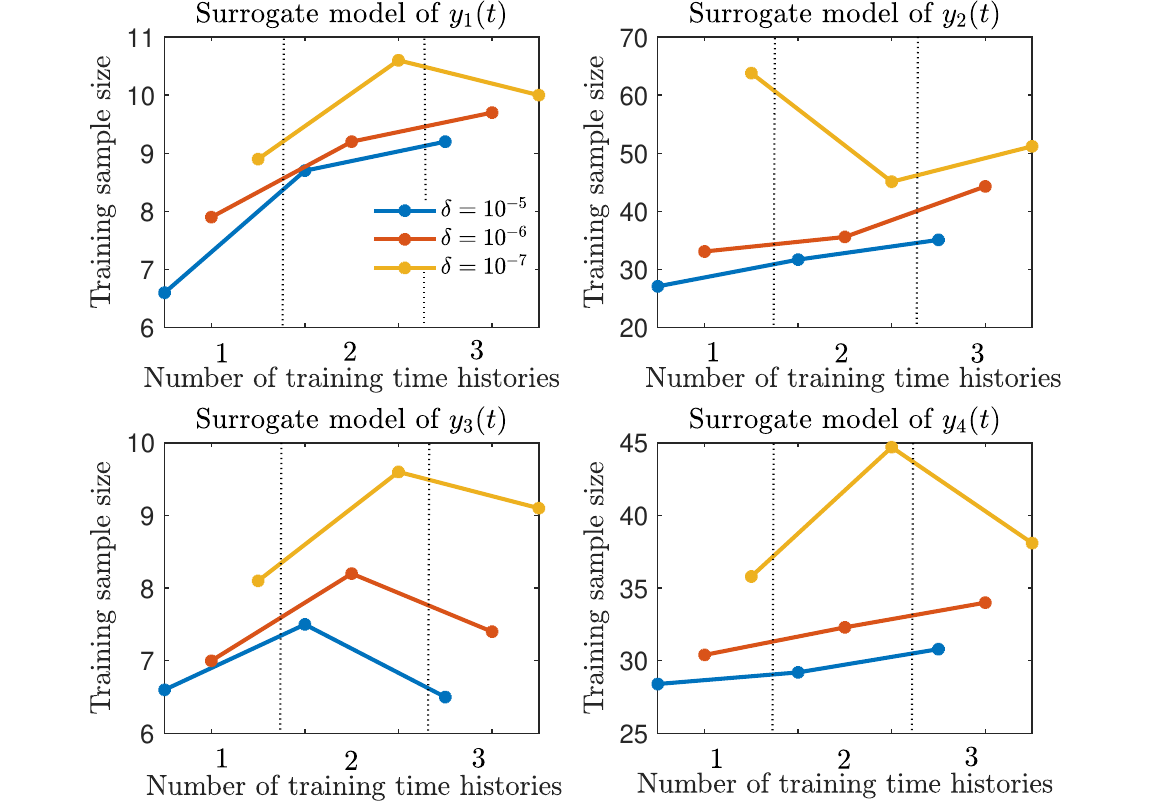}
\caption{ Sample size (mean value of 10 repetitions) of the four Kriging models in the S2K model for the quarter car model by varying the number of training time histories with different convergence thresholds.}  
\label{fig:qurater_car_sample}
\end{figure}

Fig. \ref{fig:qurater_car_error} presents the average relative error $\bar\epsilon_1$ of $x_1(t)$ predicted by the S2K model with varying the number of time histories and various convergence thresholds. It shows that the S2K model with only one training time history of the state vector already yields very accurate results. As expected, the accuracy can be further improved by increasing the number of the training time histories as well as by decreasing the convergence threshold. In addition, the training sample sizes of the four Kriging models corresponding to the four components of the state space representation in Eq. \eqref{eq:quarter_car_state_sapce} are shown in Fig. \ref{fig:qurater_car_sample}. It is observed that only a limited number of samples are retained in the four Kriging models for various settings, which confirms the sparsity of the S2K model. 

\begin{figure}[htbp]
\centering
\includegraphics[width=1\textwidth,trim={20 275 20 280},clip]
{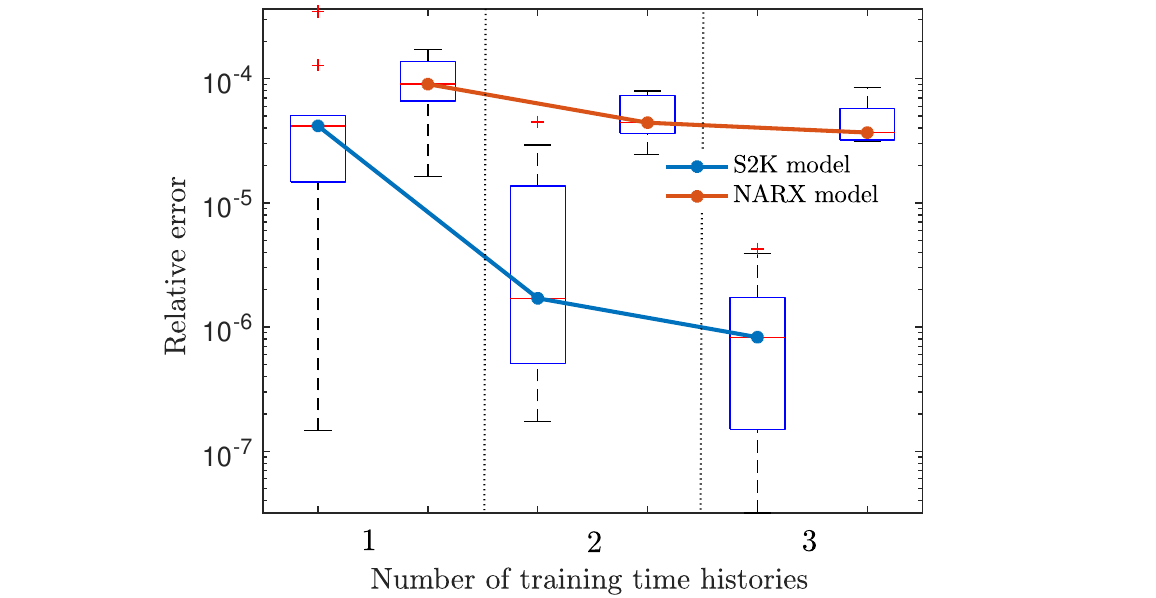}
\caption{ Comparison of the relative errors (with 10 repetitions) of $x_1$ obtained with the S2K model ($\delta=10^{-5}$) and the NARX model for the quarter car model by varying the number of training time histories. }   
\label{fig:qurater_car_comparison}
\end{figure}

The comparisons of relative error between the S2K ($\delta=10^{-5}$) and the NARX model (maximum polynomial order set to 3 \cite{mai2016surrogate}) is shown 
in Fig. \ref{fig:qurater_car_comparison}. It shows that both methods provide very accurate prediction, but the S2K outperforms the NARX model, especially when more training time histories are available.

\subsection{ Duffing oscillator}

The second example is a Duffing oscillator \cite{kougioumtzoglou2009approximate,mai2016surrogate} subjected to random loading for a duration of $T=10s$. The oscillator is governed by the following ODE
\begin{eqnarray}
\begin{aligned}
 \Ddot{x}(t) + 2\zeta\omega_n\dot x(t) + \omega^2_nx(t) +\beta x^3(t) = u(t),   \\
\end{aligned}
\label{eq:duffing}
\end{eqnarray}
where $\zeta=0.05$,  $\omega_n=10$ rad/s, and $\beta=2000 \enspace\mathrm{m^{-2}s^{-2}}$. 
The random loading $u(t)$ is modeled by a white noise ground acceleration discretized
in the frequency domain as \cite{shinozuka1991simulation}
\begin{eqnarray}
\begin{aligned}
 u(t) = \sqrt{2S\triangle\omega}\sum_{i=1}^{d/2} [\vartheta_i \cos(\omega_it) + \vartheta_{d/2+i} \sin(\omega_it) ],   \\
\end{aligned}
\label{eq:excitation_frequency}
\end{eqnarray}
where $\omega_i=i \triangle\omega$ with $\triangle\omega=30\pi/d$, and $\vartheta_i\sim \mathcal {N}(0,1)$, for $i=1,...,d$. In the current work, the spectral intensity of the Gaussian white noise is set to $S=0.1 \enspace\mathrm{m^{2}/s^3}$, and we choose $d=150$ \cite{PAPAIOANNOU2019106564} . 

Denoting  $\bsx(t)=[x(t),\dot x(t)]^{\rm T}$ and $\bsy(t)=\dot\bsx(t)$, the original differential equation in Eq. \eqref{eq:duffing} can be expressed by its state space representation as
\begin{eqnarray}
\begin{aligned}
\begin{cases} 
{y}_1(t) &  =\dot{x}(t), \\  {y}_2(t)  &=  u(t) -  2\zeta\omega_n\dot x(t) - \omega^2_nx(t) -\beta x^3(t). 
\end{cases}
\label{eq:state_duffing}
\end{aligned}
\end{eqnarray}

\begin{figure}[htp]
\centering
\includegraphics[width=0.9\textwidth,trim={1 155 20 150},clip]{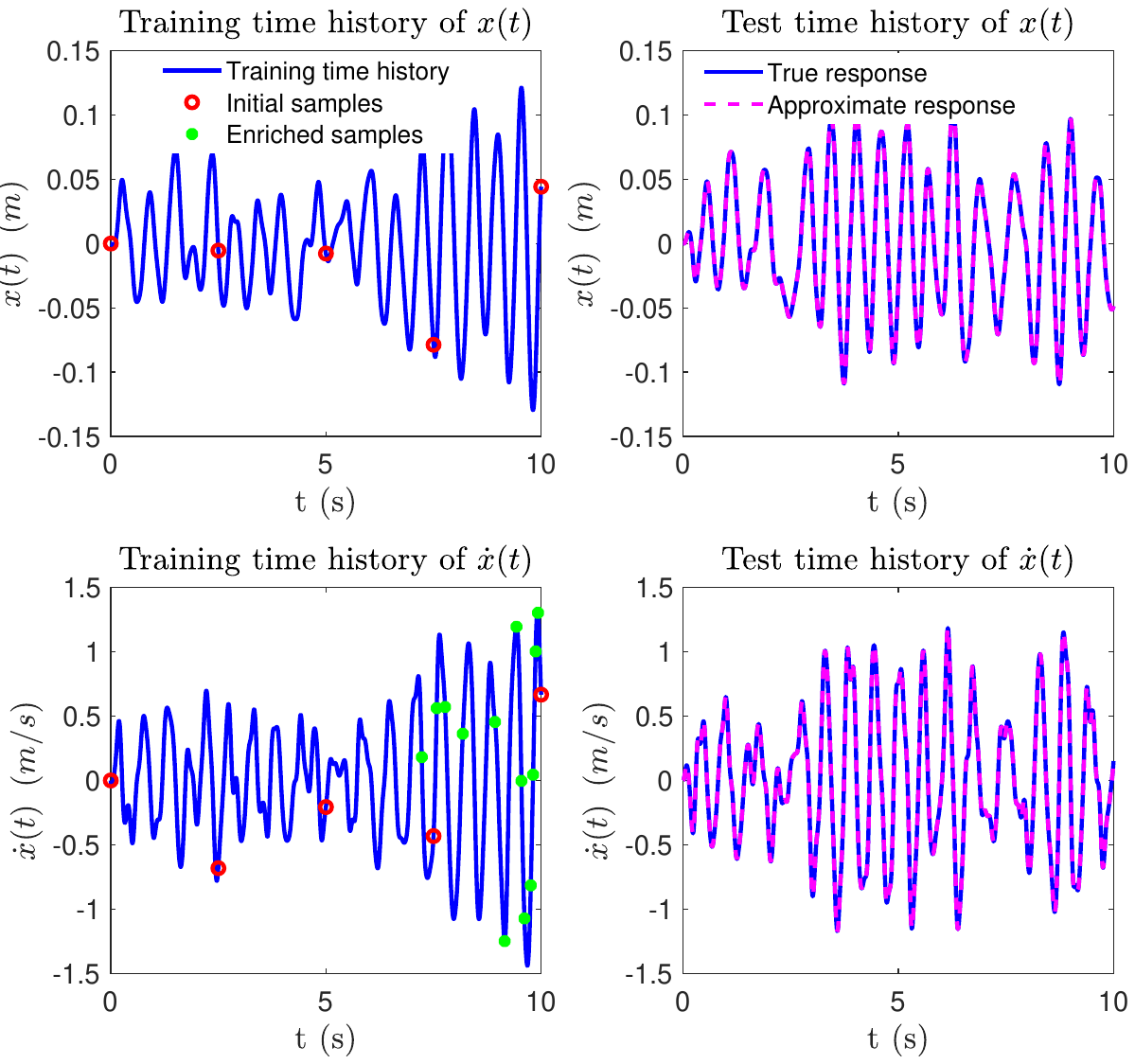}
\caption{ An example of the training and test time history of the state vector $\bsx(t)=[x(t),\dot x(t)]^{\rm T}$ of the Duffing oscillator.}
\label{fig:duffing_example}
\end{figure}

The S2K model is used to learn the dynamics of the Duffing oscillator in Eq. \eqref{eq:state_duffing}. In this example, we use only one time history of the state vector estimated for the total duration $T = 10$ s with time step $\triangle t = 0.002$ s to train the two Kriging models separately.
The results are depicted in Fig. \ref{fig:duffing_example}, where a test time history is also presented. One can see that the S2K model is highly accurate in emulating the Duffing model, with relative errors ${\epsilon}_1=9.68\times 10^{-6}$ for $x(t)$ and ${\epsilon}_2=9.09\times 10^{-6}$ for $\dot x(t)$. Again, only a few informative samples are enriched during the active learning procedure, which leads to highly sparse Kriging models. 

\begin{figure}[htp]
\centering
\includegraphics[width=1\textwidth,trim={20 255 20 255},clip]{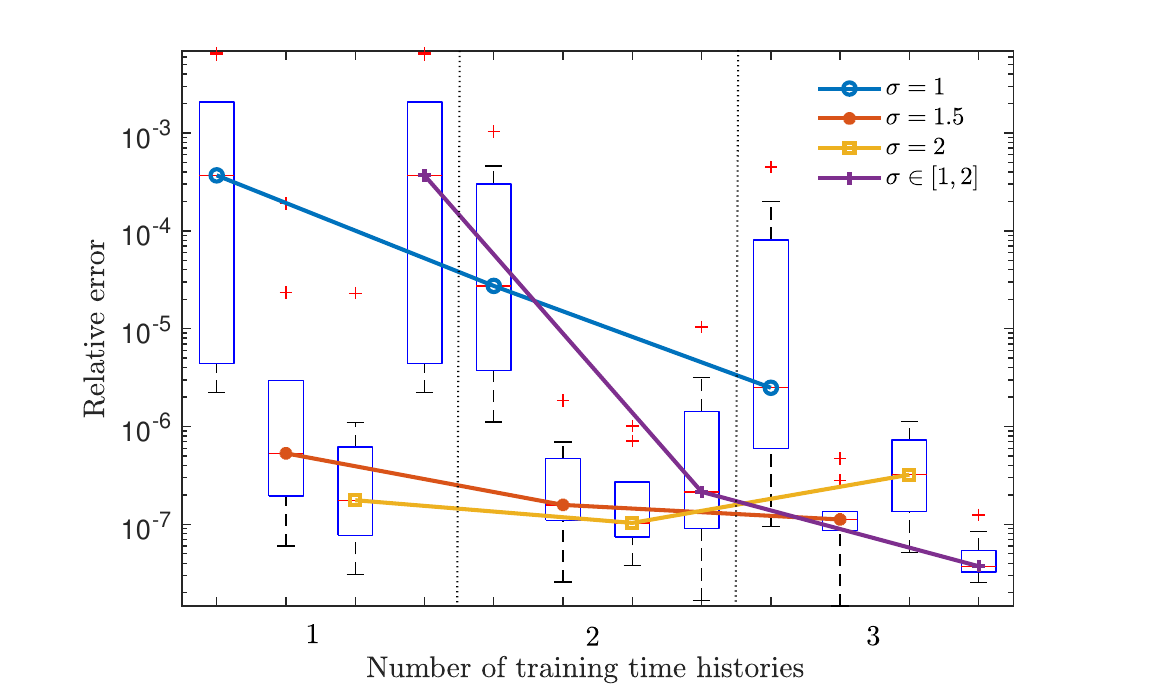}
\caption{ Boxplot of relative errors (with 10 repetitions) of $x(t)$ of Duffing oscillator obtained by the S2K model ($\delta=10^{-6}$) by varying the  number of training time histories with different magnification factors. }
\label{fig:duffing_error}
\end{figure}

The average relative error $\bar\epsilon_1$ of $x(t)$ predicted by the S2K model by varying the training time histories with different variability magnification factors of excitation are presented in Fig. \ref{fig:duffing_error}, in which $\sigma\in [1,2]$ signifies that different magnification factors are assigned to different training time histories according to Eq. \eqref{eq:mixture_sigma}. The variability of the excitation is magnified by amplifying the standard deviation of $\vartheta_i\enspace(i=1,...,150)$ in Eq. \eqref{eq:excitation_frequency}.
The S2K model with only one training time history of the state vector already provides prediction  with less than $1\%$ relative error, and its accuracy can be improved by about two orders of magnitude by magnifying the variability of the excitation by $\sigma=1.5$ or $\sigma=2$. The mixture of $3$ training time histories of state vector with different $\sigma$ values yields the optimal results for $n=3$.

\begin{figure}[htp]
\centering
\includegraphics[width=1\textwidth,trim={0 280 5 290},clip]{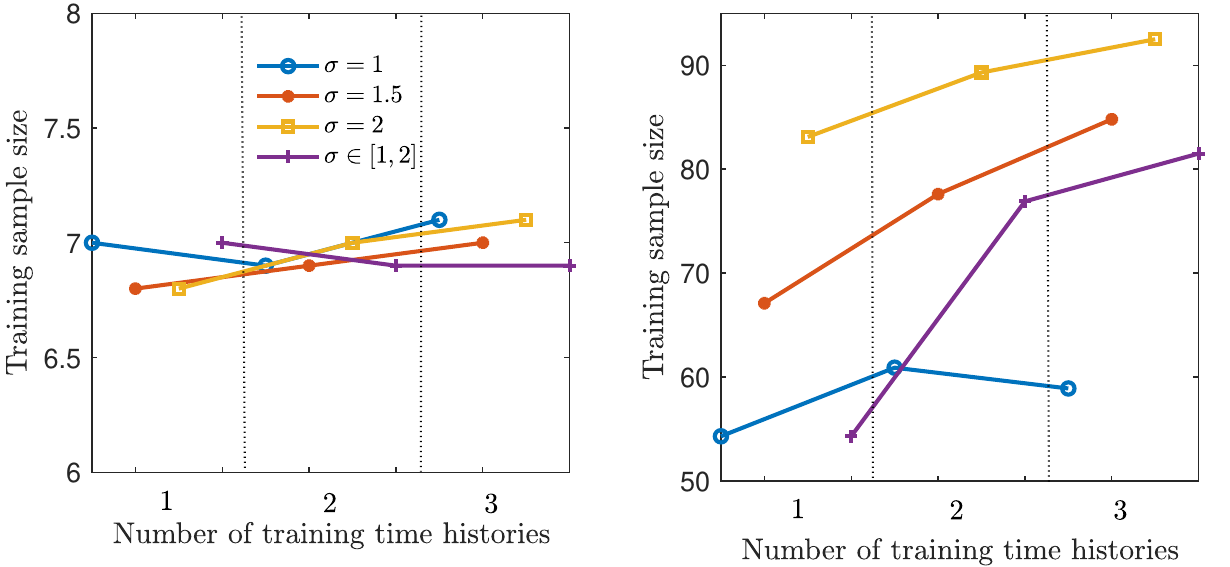}
\caption{ Training sample size (mean value of 10 repetitions) of the two Kriging models (left: $\hat y_1(t)$, right: $\hat y_2(t)$) in the S2K model ($\delta=10^{-6}$) for the Duffing oscillator by varying the number of training time histories with different magnification factors.}
\label{fig:duffing_sample}
\end{figure}

The training sample sizes of the two Kriging models corresponding to the two components of the state space representation in Eq. \eqref{eq:state_duffing} after the active learning procedure are depicted in Fig. \ref{fig:duffing_sample}. As expected, more samples are selected to train the Kriging model when using a larger magnification factor $\sigma$. However, both Kriging models remain  sparse for all settings. 

\begin{figure}[htp]
\centering
\includegraphics[width=1\textwidth,trim={20 280 20 280},clip]{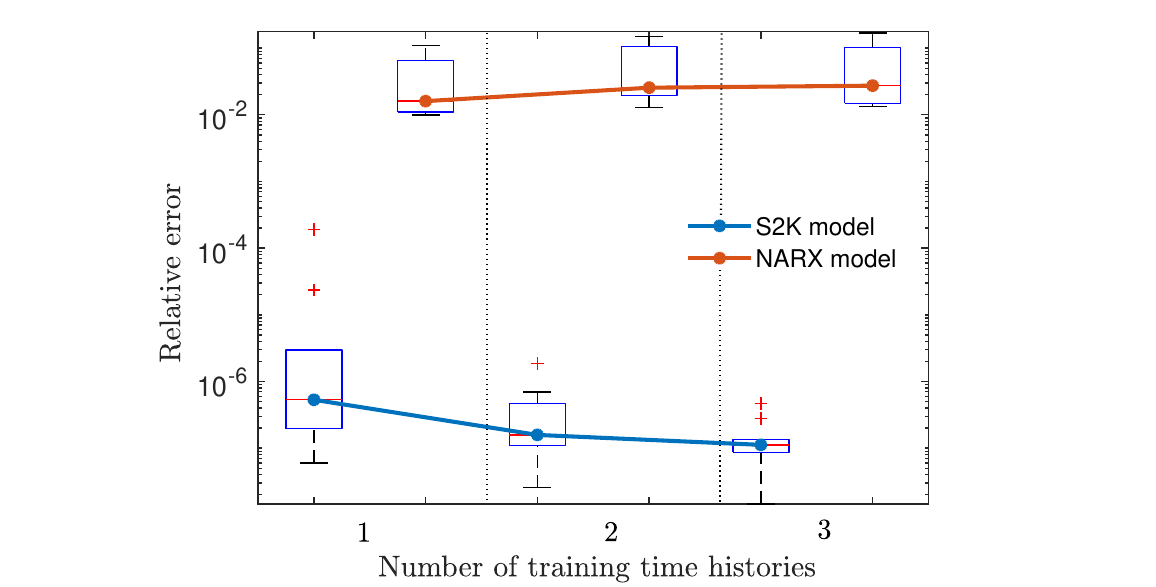}
\caption{ Comparison of the relative errors (with 10 repetitions) obtained with the S2K model ($\delta=10^{-6},\sigma=1.5$) and the NARX model for Duffing oscillator by varying the number of training time histories. }
\label{fig:duffing_comparison}
\end{figure}

In Fig. \ref{fig:duffing_comparison}, we compare the relative errors of $x(t)$ obtained with the S2K model ($\delta=10^{-6},\sigma=1.5$) and the NARX model (maximum polynomial order set as 3 \cite{mai2016surrogate}) by varying the number of training histories. Due to the strong variability of the state quantity $x(t)$ over the time domain, the NARX model gives poor prediction.

\subsection{Nonlinear Bouc-Wen hysteretic oscillator}
\label{subsec:bouc}
In the third example, we consider a Bouc-Wen hysteretic oscillator under random external excitation \cite{kafali2007seismic,PAPAIOANNOU2019106564,mai2016surrogate}, described by the following differential equation
\begin{eqnarray}
\begin{aligned}
\begin{cases} 
m \Ddot{x}(t) + c\dot x(t) + k\left [\alpha x(t) +(1-\alpha)x_y z(t)\right] = m u(t), \\
 \\  \dot z(t) = \frac{1}{x_y}\left [A\dot x(t) -\beta |\dot x(t)|| z(t)|^{d-1}z(t) -\gamma \dot x(t) |z(t)|^d \right],
\end{cases}
\end{aligned}
\label{eq:bouc_wen}
\end{eqnarray}
where the mass, stiffness and damping of the oscillator are 
$m=6\times10^4\enspace\mathrm{kg}$, $k=5\times10^6\enspace\mathrm{N}$ and $c = 2m\zeta\sqrt{k/m}$, respectively, with $\zeta=0.05$. The degree of
hysteresis is defined by $\alpha$, which is chosen as $0.5$. In addition, we set $x_y=0.04\enspace\mathrm{m}$, $\beta =\gamma = 0.5$ and $A=1,d=3$.

Denoting  $\bsx(t)=[x(t),\dot x(t),z(t)]^{\rm T}$ and $\bsy(t)=\dot\bsx(t)$, the above differential equation of Bouc-Wen oscillator can be expressed by its state space form as 
\begin{eqnarray*}
\begin{aligned}
\begin{cases} 
{y}_1(t) &  =\dot{x}(t), \\  {y}_2(t)  &= -\frac{c}{m}\dot x(t) -\frac{k}{m}\left [\alpha x(t) +(1-\alpha)x_y z(t)\right] + u(t),  \\
{y}_3(t)  &= \frac{1}{x_y}\left [A\dot x(t) -\beta |\dot x(t)|| z(t)|^{n-1}z(t) -\gamma \dot x(t) |z(t)|^n \right].
\label{eq:bouc_state}
\end{cases}
\end{aligned}
\end{eqnarray*}

The excitation is modeled by a white noise ground acceleration discretized in frequency domain as in Eq. \eqref{eq:excitation_frequency} with spectral intensity being $S=0.05 \enspace\mathrm{m^2/s^3}$. The time period is set to $T=8$ s \cite{PAPAIOANNOU2019106564}. 

\begin{figure}[htp]
\centering
\includegraphics[width=0.85\textwidth,trim={1 278 1 280},clip]{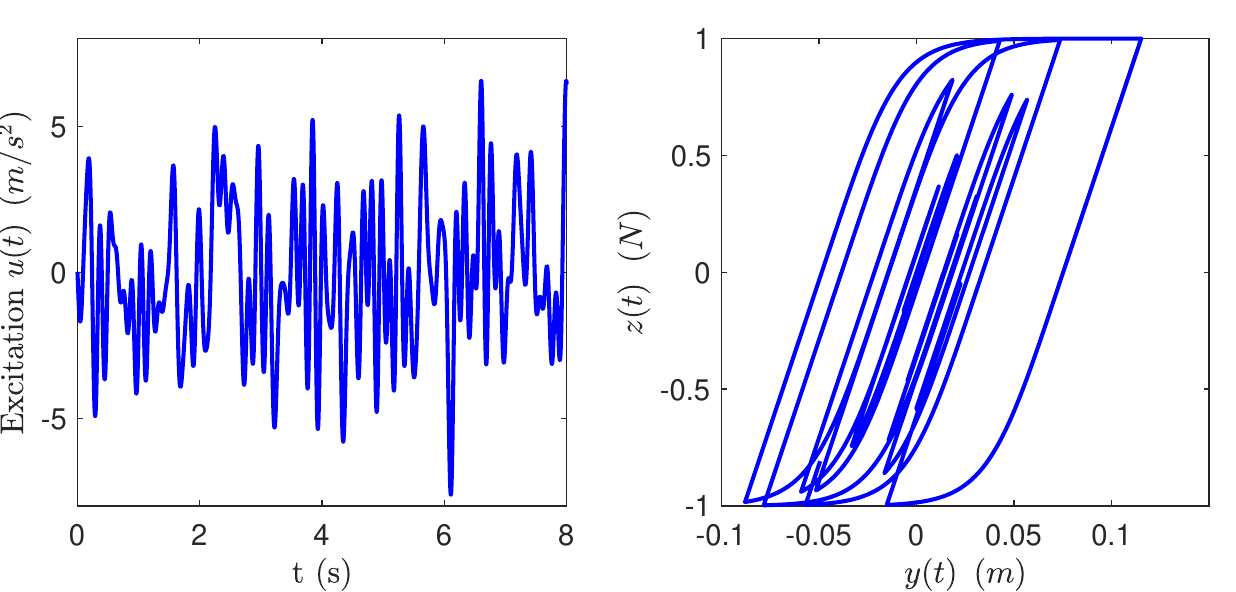}
\caption{ An example of the excitation of Bouc-Wen model }
\label{fig:Bouc_Wen_nonlinear}
\end{figure}

\begin{figure}[htp]
\centering
\includegraphics[width=1\textwidth,trim={30 155 30 155},clip]{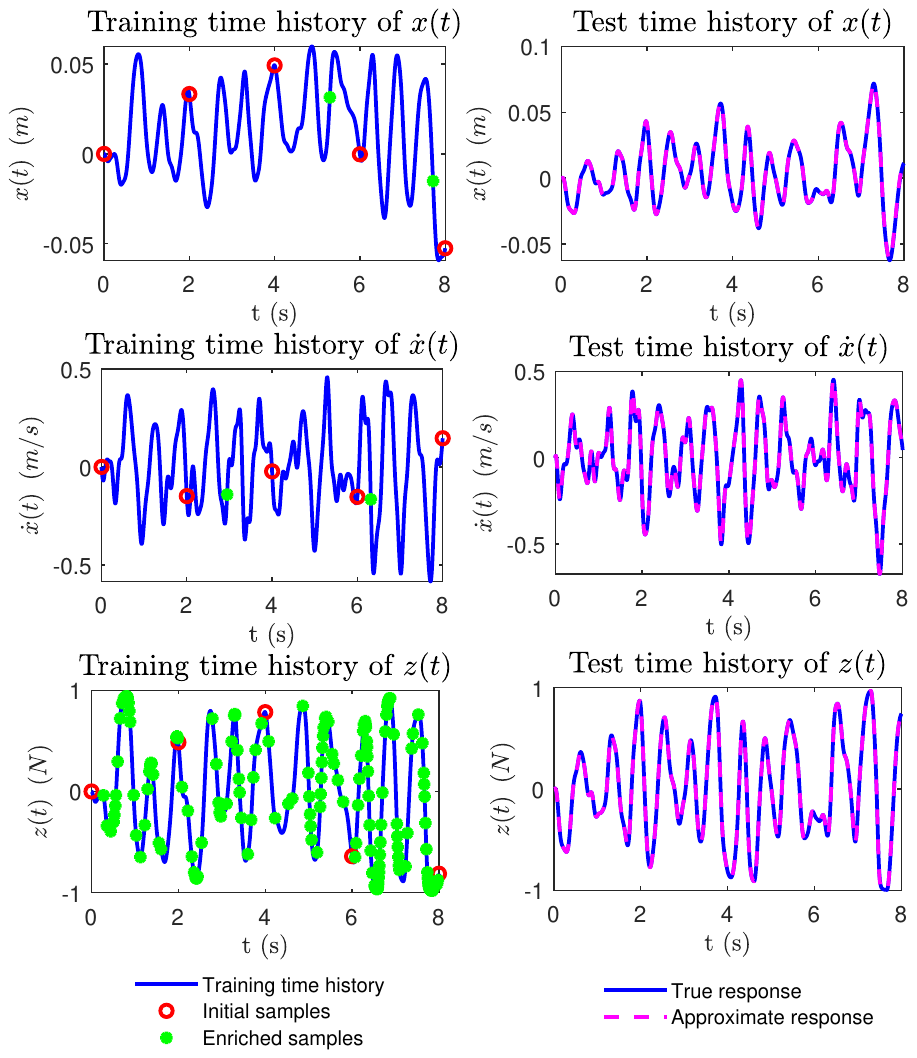}
\caption{ An example of the training and test time history of the state vector $\bsx(t)=[x(t),\dot x(t),z(t)]^{\rm T}$ of the Bouc-Wen model }
\label{fig:Bouc_wen}
\end{figure}

We first draw one time history of excitation, and the state vector is estimated with  time step $\triangle t=0.002$ s. The corresponding excitation is depicted in Fig. \ref{fig:Bouc_Wen_nonlinear}, together with the nonlinear response of the auxiliary variable $z$ of the Bouc-Wen model. We use this time history to train the S2K model with  convergence threshold $\delta=10^{-4}$. The results are depicted in Fig. \ref{fig:Bouc_wen}, in which a test time history of state vector is also presented. Again, it is shown that the S2K model yields accurate prediction, with relative errors corresponding to the test time history of state vector depicted in Fig. \ref{fig:Bouc_wen} of  ${\epsilon}_1=1.44\times 10^{-4}$ for $x(t)$, ${\epsilon}_2=1.65\times 10^{-5}$ for $\dot x(t)$ and ${\epsilon}_3=2.47\times 10^{-5}$ for $z(t)$. Since the third state equation in Eq. \eqref{eq:bouc_state} is highly nonlinear, 173 samples (5 initial samples plus 168 enriched samples) are selected from the whole time history of $z(t)$ to construct the third Kriging model. 

\begin{figure}[htp]
\centering
\includegraphics[width=1\textwidth,trim={10 265 10 265},clip]{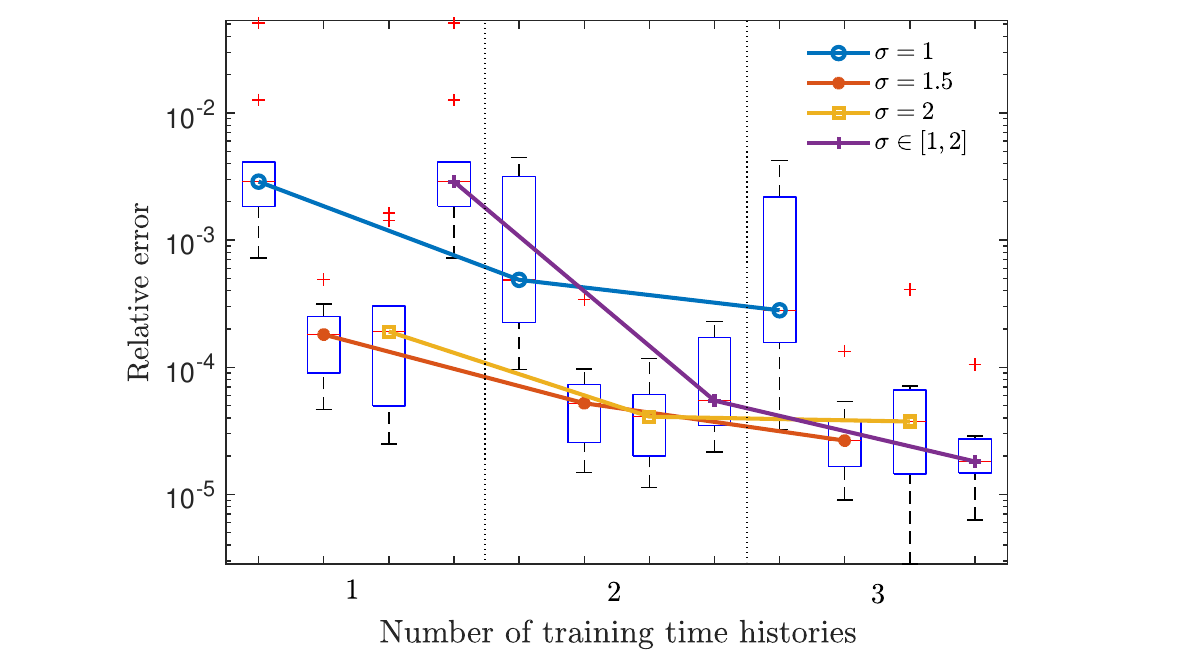}
\caption{ Boxplot of relative errors of $x(t)$  (with 10 repetitions) of the Bouc-Wen model obtained with the S2K model ($\delta=10^{-4}$) by varying the number of training time histories with different magnification factors $\sigma$. }
\label{fig:bouc_error}
\end{figure}

The average relative error $\bar\epsilon_1$ of $x(t)$ predicted by the S2K model by varying the training time histories with different 
magnification factors is presented in Fig. \ref{fig:bouc_error}, in which $\sigma \in [1, 2]$  signifies that different magnification factors are assigned to different training time histories
according to Eq. \eqref{eq:mixture_sigma}. One can see that the S2K model with only one training time history of the state vector already yields accurate predictions. Its accuracy can be improved by about one order of magnitude by magnifying the variability of the excitation by $\sigma=1.5$ times, but there is no additional improvement if a larger magnification factor $\sigma=2$ is used. In addition, the accuracy can be further improved by collecting more time histories of the state quantities. Again, one can see that the mixture of the $3$ training time histories of state vector with different $\sigma$ values yields the most accurate S2K model when $n=3$.

\begin{figure}[htp]
\centering
\includegraphics[width=1\textwidth,trim={1 325 1 325},clip]{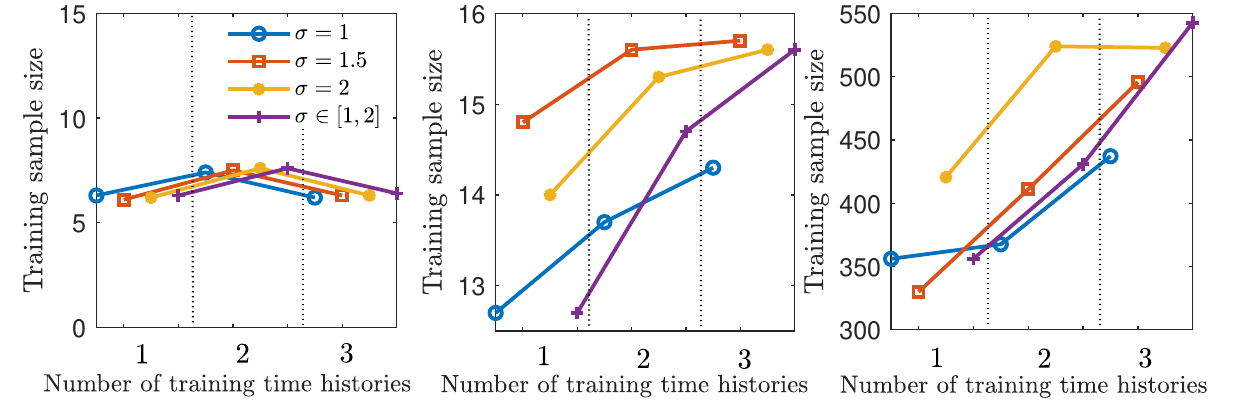}
\caption{ Sample size (mean value of 10 repetitions) of the three Kriging models (left: $\hat y_1(t)$, middle: $\hat y_2(t)$, right: $\hat y_3(t)$)  in the S2K model ($\delta=10^{-4}$) for Bouc-Wen model by varying the number of training time histories for different magnification factors $\sigma$. }
\label{fig:bouc_sample}
\end{figure}

The training sample size of the three Kriging models corresponding to the three components of the state space representation in Eq. \eqref{eq:bouc_state} after the active learning process is depicted in Fig. \ref{fig:bouc_sample}. It shows the first two Kriging models remain sparse for various parameter settings while more samples are selected to train the third Kriging model.

\begin{figure}[htp]
\centering
\includegraphics[width=1\textwidth,trim={10 280 10 285},clip]{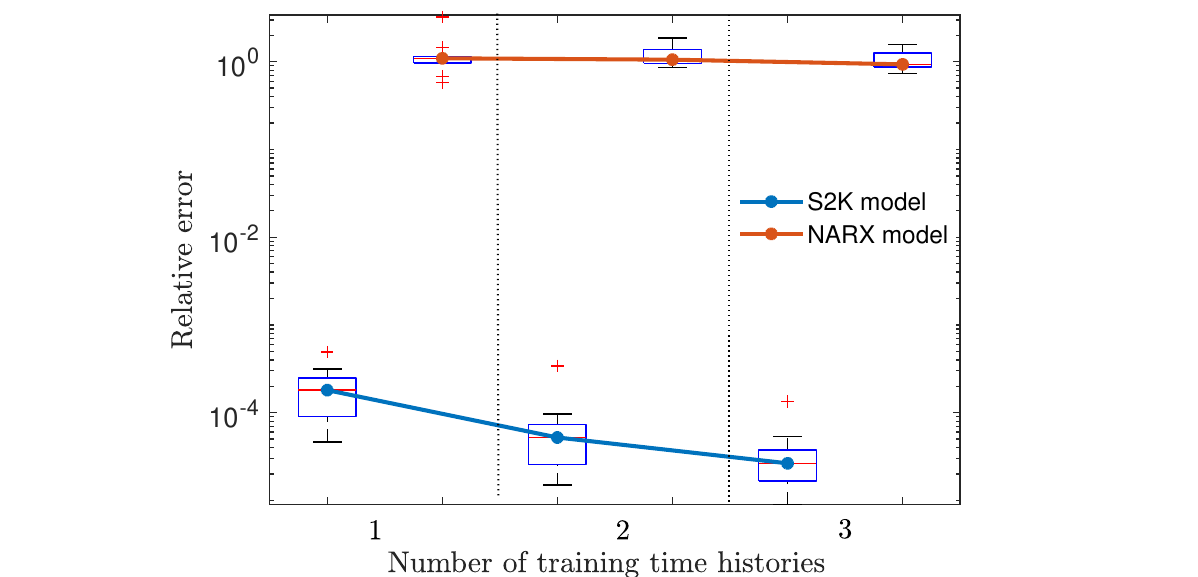}
\caption{ Comparison of relative errors (with 10 repetitions) of $x(t)$ of Bouc-Wen model obtained with the S2K model ($\delta=10^{-4},\sigma=1.5$) 
 and the NARX model by varying the number of training time histories.}
\label{fig:bouc_comparison}
\end{figure}

The comparison of the relative errors of $x(t)$ obtained with both the S2K model ($\delta=10^{-4},\sigma=1.5$) and the NARX model (maximum polynomial order 5) is depicted in Fig. \ref{fig:bouc_comparison}.
Note that different from \cite{mai2016surrogate}, we do not use any physical information to construct the NARX model here. For this challenging nonlinear problem, the vanilla NARX model with polynomial basis fails to provide satisfactory prediction. By contrast, the S2K model yields accurate prediction with relative errors that are 4 orders of magnitude lower than that of the NARX model.

\subsection{Two-story nonlinear hysteretic structure }
\label{subsec:two-story}

 We apply the S2K model to emulate the behavior of a two-story two-span nonlinear hysteretic frame structure subjected to ground motion excitation \cite{lyu2022unified}. The corresponding motion equation reads
\begin{eqnarray}
\begin{aligned}
\boldsymbol{M} \Ddot{\bsx}(t) + \boldsymbol{C} \dot{\bsx}(t) +\boldsymbol{G}[\bsx(t),\bsz(t)] = -\boldsymbol{M} \boldsymbol{I} u(t),
\end{aligned}
\label{eq:two-story}
\end{eqnarray}
where $\bsx(t)\in \mathbbm{R}^{2},\dot{\bsx}(t)\in \mathbbm{R}^{2}$, and $\Ddot{\bsx}(t)\in \mathbbm{R}^{2}$ represent the displacement, velocity and acceleration vector; $\boldsymbol{M}\in \mathbbm{R}^{2\times2}$,  $\boldsymbol{C} \in \mathbbm{R}^{2\times2}$ and $\boldsymbol{I}\in \mathbbm{R}^{2\times2}$ are the lumped mass, Rayleigh damping and the identity matrices; $u(t)$ is modeled by a  white noise ground acceleration process as that in Eq. \eqref{eq:excitation_frequency}; $\boldsymbol{G}(\cdot)$ is the restoring force vector characterized by the Bouc-Wen model as 
\begin{eqnarray*}
\tilde {G}_j(\bsx(t),\bsz(t)) = \alpha K_j \tilde {x}_j(t) + (1-\alpha)K_j z_j(t),
\end{eqnarray*}
in which $\alpha = 0.04$; $\tilde {G}_j(\cdot)$, $\tilde {x}_j(t)$, and $K_j$ represent the $j$-th inter-story restoring force, drift, and story initial stiffness, respectively; $z_j(t)(j = 1,2)$ is the $j$-th auxiliary hysteretic displacement, which is described by
\begin{eqnarray}
 \dot z_j(t) = \frac{\varpi}{1+d_{\eta}\epsilon_j(t)}\left(1-\xi {\rm exp}\left(-\left\{{\frac{z_j(t){\rm sgn}[\dot{\tilde {x}}_j(t)]-\frac{q}{(\beta+\gamma)\kappa } }{[\psi + d_{\psi}\epsilon_j(t) ](\lambda +\zeta_s\xi )}}\right\}^2\right)\right)
\label{eq:hysteretic_displacement},
\end{eqnarray}
where $\kappa =1+d_{\nu}\epsilon_j(t)$, $\varpi = \dot{\tilde {x}}_j(t) - \kappa\left[\beta |\dot{\tilde {x}}_j(t)|z_j(t) + \gamma \dot{\tilde {x}}_j(t) |z_j(t)| \right]$, $\xi= \zeta_s[1-e^{-p\epsilon_j(t)}]$; the parameters $\beta$ and $\gamma$ control the basic hysteresis shape; $d_{\nu}$ and $d_{\eta}$ are the strength and stiffness degradation; $\zeta_s$ measures the total slip; $q,p,\psi,d_{\psi}$, and $\lambda$ are the initiation, slope, magnitude, rate, and severity interaction of pinching \cite{ma2004parameter}; 
$\epsilon_j(t)$ is the $j$-th story hysteretic-dissipated energy, which is described by
\begin{eqnarray*}
 \dot\epsilon_j(t) = \dot{\tilde {x}}_j(t)z_j(t), \mathrm{for }\,\ j=1,2.
\end{eqnarray*}

In the current work, the system parameters are set as follows: lumped mass $M_j=2.6\times10^5 \,\ \mathrm{kg}$, initial stiffness $ K_j = 10^8 \,\ \mathrm{N/m}$, $\beta=15\,\ \mathrm{m^{-1}}$, $\gamma=150\,\ \mathrm{m^{-1}}$,  $d_{\nu}=d_{\eta}=p=1000\,\ \mathrm{m^{-2}}$, $q=0.25$, $d_{\psi}= 5 \,\ \mathrm{m^{-2}}$, $\lambda=0.5$, $\zeta_s=0.99$, $\psi=0.05\,\ \mathrm{m}$. In addition, the damping ratios of the first two modes are taken as $\zeta_1=\zeta_2=0.05$.
 
\begin{figure}[htp]
\centering
\includegraphics[width=1.05\textwidth,trim={1 275 1 270},clip]{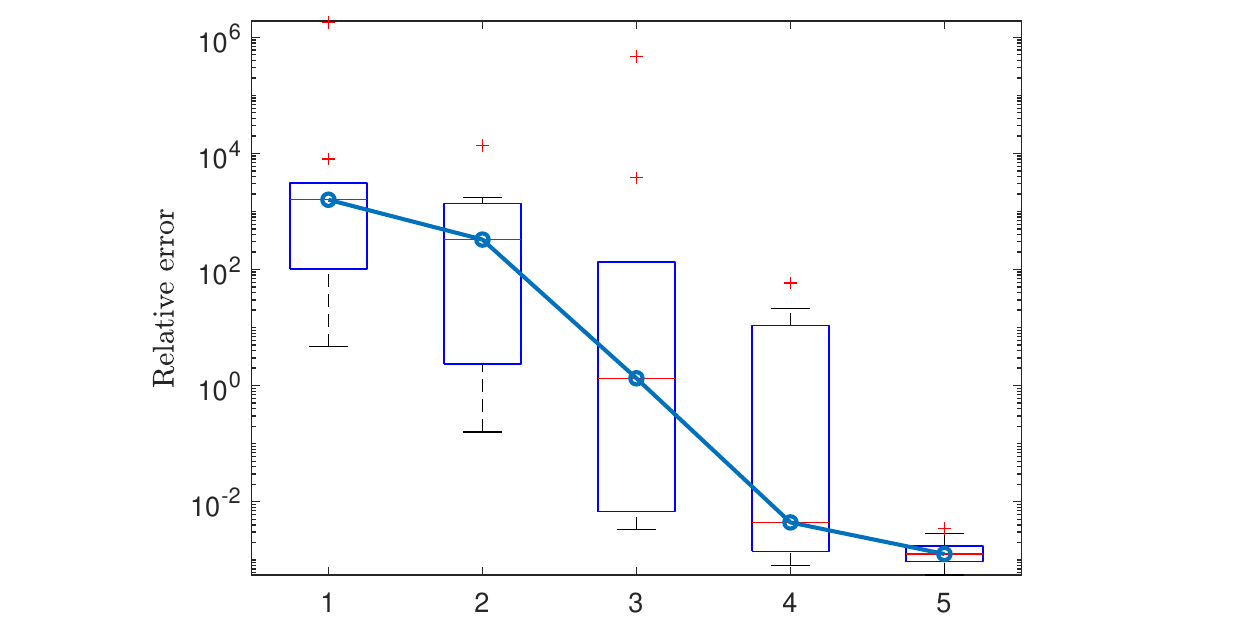}
\caption{ Boxplot of relative error in the S2K model ($\delta=5\times10^{-4}$) for the two-story nonlinear hysteretic structure by varying the number of training time histories in Example \ref{subsec:two-story}. }
\label{fig:2dof_error}
\end{figure}

\begin{figure}[htp]
\centering
\includegraphics[width=1.05\textwidth,trim={40 195 30 168},clip]{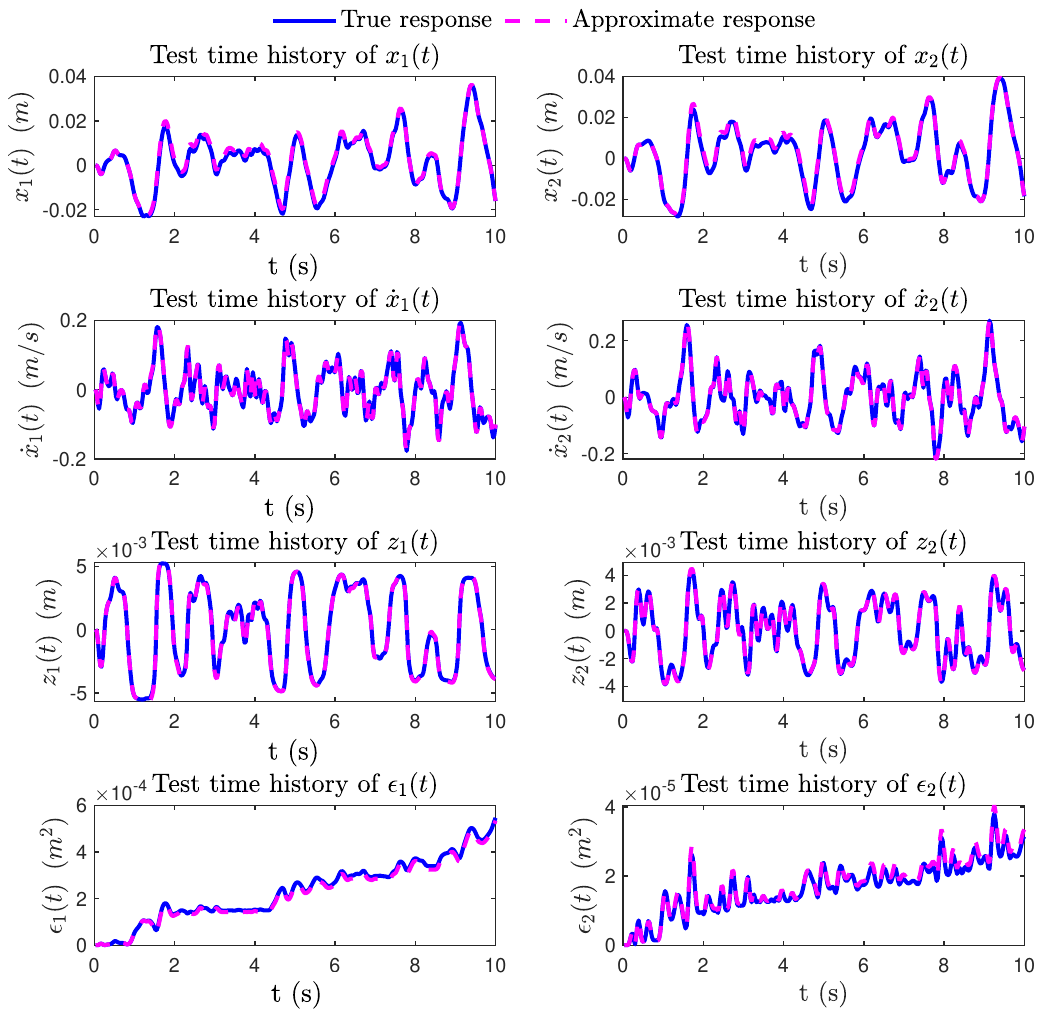}
\caption{ A test time history of the 8-dimensional state vector  $\bsx(t)=[x_1(t),x_2(t),\dot x_1(t),\dot x_2(t),z_1(t), z_2(t), \epsilon_1(t), \epsilon_2(t)]^{\rm T}$ in Example \ref{subsec:two-story}. }
\label{fig:2dof_test}
\end{figure}

\begin{figure}[htp]
\centering
\includegraphics[width=1.05\textwidth,trim={1 240 1 240},clip]{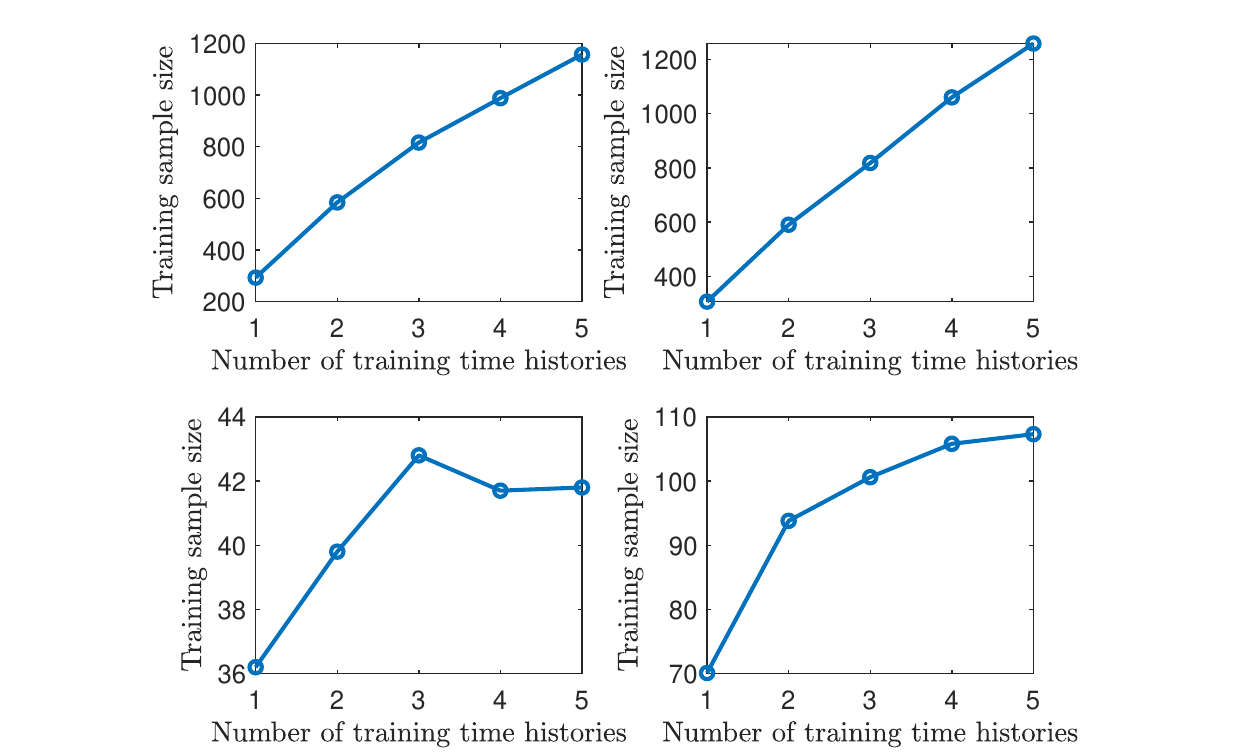}
\caption{ Sample size  (mean value of 10 repetitions) of the last four Kriging models (top left : $\hat y_5(t)$, top right : $\hat y_6(t)$, left bottom : $\hat y_7(t)$, right bottom : $\hat y_8(t)$)  in the S2K model ($\delta=5\times10^{-4}$) for the two-story nonlinear hysteretic structure by varying the number of training time histories. }
\label{fig:2dof_sample}
\end{figure}

The state vector $\bsx(t)=[x_1(t),x_2(t),\dot x_1(t),\dot x_2(t),z_1(t), z_2(t), \epsilon_1(t), \epsilon_2(t)]^{\rm T}$ is 8-dimensional, and we construct 8 Kriging models in the S2K model to learn the state space representation of this dynamical system, where the input is $[\bsx(t)^{\rm T}, u(t)]^{\rm T}$, and the output is $\bsy(t)=\dot\bsx(t)$. The time history of the state vector over the whole time period of interest $T\in [0,10]$ s is evaluated with time step $\triangle t = 0.01$ s. The average relative errors $\bar\epsilon_1$ of $x_1(t)$ predicted by the S2K model with varying number of training time histories are presented in Fig. \ref{fig:2dof_error}, in which  different magnification factors are assigned to different training time histories
according to Eq. \eqref{eq:mixture_sigma}. Due to the strong non-linearity and non-smoothness property of this model, at least 5 time histories of the state vector are required to train an accurate S2K model. To illustrate its performance, a specific test time history of the state vector is depicted in Fig. \ref{fig:2dof_test}. It shows that the eight state quantities are well predicted by the S2K model over the whole time period of interest (with 5 training time histories).  In this example, only the last four components of the state space representation are strongly nonlinear, especially  $y_6(t)=\dot z_1(t)$ and $y_7(t)=\dot z_2(t)$ in Eq. \eqref{eq:hysteretic_displacement}.  We therefore only present the training sample sizes of the last four Kriging models in Fig. \ref{fig:2dof_sample}. As expected, as the number of training time histories increases, an increasing number of samples is selected during the active learning procedure to train the Kriging models, and more than $20 \%$ of the total samples from the 5 training time histories are selected to train the two challenging functions: $y_6(t)=\dot z_1(t)$ and $y_7(t)=\dot z_2(t)$ in Eq. \eqref{eq:hysteretic_displacement}. In this example, application of the NARX model led to meaningless predictions after some time period, it is not suitable for emulating the response of this system due to the high nonlinearity in Eq. \eqref{eq:hysteretic_displacement}. Hence we do not report the relative error of the NARX model for comparisons.

\section{Concluding remarks}
\label{sec:conclusion}

 In this work, we have presented a novel surrogate modeling framework for emulating dynamical systems under stochastic excitation by learning its state space representation with Kriging (S2K) model. Several conclusions can be drawn from the work:
 
 (1) Learning the state space representation of a dynamical system under stochastic excitation can avoid the curse of dimensionality of surrogate model training due to discretization of the stochastic excitation. 

 (2) The proposed active learning algorithm can select an informative sample subset from the whole training sample set efficiently, resulting in a sparse Kriging model. 
 
 (3) The proposed technique for designing the training time history of the state vector by magnifying the variability of the excitation can improve the accuracy of S2K model. 
 
 (4) Numerical examples demonstrate that the S2K model is powerful for emulating various complex nonlinear dynamical systems under stochastic excitation. It provides highly accurate prediction (Relative error is less than $10^{-3}$) with only a few time histories of the state vector, and it outperforms the NARX model in terms of accuracy and efficiency. 
 
In the future, we plan to adapt this framework to emulating complex nonlinear dynamical systems with both random system parameters and stochastic excitation. Moreover, we will also consider combining the S2K model with model reduction techniques, e.g., proper orthogonal decomposition or auto-encoder, to emulate stochastic dynamical systems with large number of degrees of freedom.

\section{Acknowledgements}
 This work was supported by the Alexander von Humboldt Foundation.

\bibliographystyle{plain}
\bibliography{Ref}
\end{document}